\input amstex
\documentstyle{amsppt}
\magnification1200
\tolerance=10000
\overfullrule=0pt
\def\n#1{\Bbb #1}
\def\p{\Bbb C_{\infty}}

\def\Ext{\hbox{Ext}}

\def\Hom{\hbox{Hom}}

\def\Tor{\hbox{Tor}}

\def\ord{\hbox{ord }}
\def\inc{\hbox{inc}}

\def\s{\sigma}

\def\e11{E_{11}}

\def\ga{\goth A}

\def\ve{\varepsilon}
\def\de{\Delta}
\def\vk{\varkappa}
\def\ga{\gamma}

\def\be{\beta}
\def\th{\theta}
\def\al{\alpha}
\def\si{\sigma}

\def\la{\lambda}
\def\vf{\varphi}

\def\g{\goth }

\topmatter
\title
$h^1\ne h_1$ for Anderson t-motives
\endtitle
\author
A. Grishkov, D. Logachev\footnotemark \footnotetext{E-mails: shuragri{\@}gmail.com; logachev94{\@}gmail.com (corresponding author)\phantom{*******************}}
\endauthor
\thanks Thanks: The authors are grateful to FAPESP, S\~ao Paulo, Brazil for a financial support (process No. 2017/19777-6). The first author is grateful to SNPq, Brazil, to RFBR, Russia, grant 16-01-00577a (Secs. 1-4), and to Russian Science Foundation, project 16-11-10002 (Secs. 5-8) for a financial support. The authors are grateful to (an) anonymous reviewer(s) for some important remarks, particularly, for the proofs of Lemma 1.4.2, affirmation of Remark 1.7.5 and its proof, indication of existence of $\omega_1$ in (1.7.5.1), and indication of references [HJ], [M], [P] where propositions similar to the ones of the present paper are proved.
\endthanks
\NoRunningHeads
\address
First author: Departamento de Matem\'atica e estatistica
Universidade de S\~ao Paulo. Rua de Mat\~ao 1010, CEP 05508-090, S\~ao Paulo, Brasil, and Omsk State University n.a. F.M.Dostoevskii. Pr. Mira 55-A, Omsk 644077, Russia.
\medskip
Second author: Departamento de Matem\'atica, Universidade Federal do Amazonas, Manaus, Brasil
\endaddress
\abstract Let $M$ be an Anderson t-motive of dimension $n$ and rank $r$. Associated are two $\n F_q[T]$-modules $H^1(M)$, $H_1(M)$ of dimensions $h^1(M)$, $h_1(M)\le r$ - analogs of $H^1(A,\n Z)$, $H_1(A,\n Z)$ for an abelian variety $A$. There is a theorem (Anderson): $h^1(M)=r \iff h_1(M)=r$; in this case $M$ is called uniformizable. It is natural to expect that always $h^1(M)=h_1(M)$. Nevertheless, we explicitly construct a counterexample. Further, we answer a question of D.Goss: is it possible that two Anderson t-motives that differ only by a nilpotent operator $N$ are of different uniformizability type, i.e. one of them is uniformizable and other not? We give an explicit example that this is possible.
\endabstract
\endtopmatter
\document
{\bf 0. Statement of the problem.}
\medskip
Let $q$ be a power of a prime number $p$ and $\n F_q$ the finite field of order $q$. The field $\n F_q(\theta)$ is its field of rational functions, and it is the function field analog of $\n Q$. There is a valuation ord on $\n F_q(\theta)$ defined by the condition $\ord \th=-1$.
\medskip
The field of the Laurent series $\n F_q((1/\theta))$ is the completion of $\n F_q(\theta)$ in the topology defined by ord. It is the function field analog of $\n R$. Let $\p$ be the completion of the algebraic closure of $\n F_q((1/\theta))$, which is the function field analog of $\n C$. By definition, $\p$ is complete. It is also algebraically closed ([G], Proposition 2.1).
\medskip
Let $M$ be an Anderson t-motive of dimension $n$ and rank $r$ over $\p$. It is a function field analog of an abelian variety, more exactly, of an abelian variety $A$ with multiplication by an imaginary quadratic field, of dimension $r$ and of signature $(n,r-n)$. Attached to such $A$ is its lattice $H_1(A)$ which is isomorphic to $\n Z^{2r}$. The cohomology group $H^1(A)$ is the dual lattice, there is a perfect pairing between $H_1(A)$ and $H^1(A)$.
\medskip
{\bf Remark 0.1.} $n$ can be considered as an analog of dimension of $A$. Really, the above $A$ can be defined by means of a "lattice" $H_1(A)$ of dimension $r$ in $\n C^n$, see [GL09] for details.\footnotemark \footnotetext{We do not need this construction in the present paper.}
\medskip
Analogs of (co)homology groups $H_1(A)$, $H^1(A)$ can be defined for Anderson t-motives. $H^1(M)$ is defined in [G], 5.9.11.2. For the definition of $H_1(M)$ we use the following notations. There exist objects called t-modules (see [AN]; [G], Definition 5.4.5); their category is anti-equivalent to the category of Anderson t-motives ([G], Theorem 5.4.11).
\medskip
{\bf 0.2.} Let $E$ be a t-module. Following [G], the line below the proof of Theorem 5.4.11, we denote the corresponding Anderson t-motive by $M=M(E)$. Analogically, let $M$ be an Anderson t-motive. We denote the corresponding t-module by $E=E(M)$.
\medskip
For any t-module $E$ the set $H_1(E)$ is defined in [G], 5.9.11.3. By definition, $H_1(M):=H_1(E(M))$ (see (1.5.1) below for details). \footnotemark \footnotetext{Unlike the standard notations, $H^1$ is a covariant functor in $M$ and $H_1$ is contravariant, because initially they were defined for the anti-equivalent category of t-modules.} Both $H^1(M)$ and $H_1(M)$ are free $\n F_q[T]$-modules (here $T$ is an abstract variable, see (1.1)). Their dimensions are denoted by $h^1(M)$, $h_1(M)$ respectively. We have $h^1(M), \ h_1(M)\le r$.
\medskip
{\bf 0.2.1.} There exists a canonical pairing $$\pi: H_1(M)\underset{\n F_q[T]}\to{\otimes}H^1(M)\to\n F_q[T]$$ (see (1.7.3), (1.7.4)). We have
\medskip
{\bf Theorem 0.2.2.} (Anderson, [AN]; [G], 5.9.14). $h^1(M)=r \iff h_1(M)=r$. In this case $\pi$ is perfect over $\n F_q[T]$.
\medskip
The t-motives satisfying this condition are called uniformizable. If $M$ is uniformizable then $H_1(M)$ is isomorphic to the lattice $L(M)$ of $M$, where $L(M)\subset \p^n$. See [G], 5.9.11.3 for details (multiplication by $T$ in $H_1(M)$ corresponds to the linear transformation $\th I_n +N$ in $L(M)$, see (1.3.1) for the definition of $N$. Particularly, if $N=0$ then the multiplication by $T$ in $H_1(M)$ corresponds to the multiplication by $\th$ in $L(M)$).
\medskip
We have a natural question. Let $M$ be any t-motive.
\medskip
{\bf Question 0.3.} (a) Is it always the case that $h_1(M)=h^1(M)$ ?
\medskip
(b) If $h_1(M)=h^1(M)$, what is the type of $\pi$? It can be either perfect over $\n F_q[T]$, or perfect only over $\n F_q(T)$, or non-perfect.
\medskip
(c) If $\pi$ is not perfect, what is its possible rank? Can it take all possible values from 0 to min ($h^1(M), \ h_1(M)$), or not? Maybe it is always equal to min ($h^1(M), \ h_1(M)$)? Maybe it is never 0 (if $h^1(M)\ne0\ne h_1(M)$)?
\medskip
 By analogy with the number field case of abelian varieties and by Theorem 0.2.2, it is natural to expect that the answer to (a) is yes. Nevertheless, we construct a counterexample (Theorems 4.1, 4.7). Also, we give an answer to a question of D. Goss (Theorem 5.1).
\medskip
{\bf Structure of the paper.} In Section 1 we give definitions and elementary properties of Anderson t-motives and their $H^1$, $H_1$. In Section 2 we define affine equations and introduce a notion of $T$-divisible $\n F_q[[T]]$-submodule of $\p[[T]]$ --- the main computational tools to calculate $h^1$ of t-motives. We state a conjecture 2.3.6 that the length of the tail (see a footnote below (2.1.2)) of an affine equation corresponding to an Anderson t-motive of dimension $n$, is equal to $n$. In Section 3 we show how to reduce a problem of finding of $h^1(M)$ for $M$ defined in (1.10) to a solution of an affine equation. In Section 4 we apply this method to an explicitly defined t-motive to get a counterexample to (0.3a). In Section 5 we consider an example of a t-motive giving answer to a question of D. Goss. We consider in Appendix some possibilities of further research. Namely, in 6.1 - 6.4 we consider initial steps of calculations of Section 3 for some other types of t-motives, in order to confirm Conjecture 2.3.6. In (6.5) we define some Ext groups; maybe they are interesting for development of this subject. We state in (6.6) the problem of finding of all possible values of $h^1$, $h_1$ and the rank of pairing, and state a conjecture on the set of these values. We justify in (6.7) this conjecture.
\medskip
\medskip
{\bf 1. Definitions.} Let $\p[T,\tau]$ be the Anderson ring, i.e. the ring of non-commutative polynomials in two variables $T$, $\tau$ over $\p$ satisfying the following relations (here
$a \in \p$):
$$Ta=aT, \ T\tau = \tau T, \ \tau a = a^q \tau. \eqno{(1.1)}$$
Subrings of $\p[T,\tau]$ generated by $\tau$, resp. $T$ are denoted by $\p\{\tau\}$ (the ring of non-commutative polynomials in one variable), resp. $\p[T]$ (the ordinary ring of (commutative) polynomials in one variable).
\medskip
There are various versions of the definitions of t-motives. For our purposes it is sufficient to consider the following
\medskip
{\bf Definition 1.2.} ([G], 5.4.2, 5.4.18, 5.4.16). A t-motive\footnotemark \footnotetext{Terminology of Anderson; Goss calls these objects abelian t-motives.} $M$ is a left
$\p[T, \tau]$-module which is free and finitely generated as both $\p[T]$-,
$\p\{\tau\}$-modules and such that
$$ \exists \goth m= \goth m(M) \ \hbox{such that}\ (T-\theta)^ \goth m M/\tau M=0.\eqno{(1.2.1)}$$
\medskip
The dimension of $M$ over $\p\{\tau\}$ (resp. $\p[T]$) is denoted by $n$ (resp.
$r$), and these numbers are called the dimension and rank of $M$.
\medskip
{\bf Remark 1.2.2.} Let $E$ be a t-module. By definition, the dimension and rank of $M(E)$ are called by the dimension and rank of $E$ itself. t-modules of dimension 1 are called Drinfeld modules, and the t-module of dimension 1 and rank 1 (over $\p$ it is unique) is called the Carlitz module. By abuse of language, the corresponding Anderson t-motives are called the Drinfeld modules, respectively the Carlitz module (although they are t-motives, and not t-modules!).
\medskip
We can consider $M$ as a $\p[T]$-module with the action of $\tau$. Hence, sometimes for $x\in M$ we shall write $\tau(x)$ (action) instead of $\tau x$ (multiplication).
\medskip
We shall need the explicit matrix description of t-motives. First, let
$e_*=(e_1, ..., e_n)^t$ (here and below $t$ means transposition) be the column vector of elements of a basis
of $M$ over $\p\{\tau\}$. There exists a matrix $\goth A\in M_n(\p\{\tau\})$ such that

$$T e_* = \goth A e_*, \ \ \goth A = \sum_{i=0}^l \goth A_i \tau^i \hbox{ where } \goth A_i
\in M_n(\p).\eqno{(1.3)}$$
Condition (1.2.1) is equivalent to the condition
$$\goth A_0=\theta I_n + N\eqno{(1.3.1)}$$
where $N$ is a
nilpotent matrix, and the
condition \{$ \goth m(M)$ can be taken to 1\} is equivalent to the condition $N=0$.
\medskip
Second, let $f_*=(f_1, ..., f_r)^t$ be the column vector of elements of a basis
of $M$ over $\p[T]$. There exists a matrix $Q\in M_r(\p[T])$ such that $$Qf_*=\tau f_*.\eqno{(1.4.1)}$$

For all t-motives we have $Q\in GL_r(\p(T))$. This follows immediately from a more exact result:
\medskip
{\bf Lemma 1.4.2.} $\det(Q) = c(T-\th)^n$ where $c\in \p^*$.
\medskip
{\bf Proof.} Condition (1.2.1) implies that $(T-\th)^\g m$
annihilates $M/\tau M$. Hence, $M/\tau M$ is isomorphic as a $\p[T]$-module to the direct sum of finite cyclic $(T-\th)$-power torsion $\p[T]$-modules. Its dimension over $\p$ is $n$. On the other hand, (1.4.1) implies that as $\p[T]$-modules, we have
$$M/\tau M\cong M_{1\times r}(\p[T]) \ / \ M_{1\times r}(\p[T])\cdot Q.$$
This implies the lemma. $\square$
\medskip
We use some definitions of [G], 5.9.10. First, we denote by $\p\{T\}$ a subring of $\p[[T]]$ formed by series $\sum_{i=0}^\infty a_iT^i$ such that lim $a_i=0 \ ( \iff \ord a_i\to +\infty$). These series will be called "small" series; the same terminology will be applied to elements of free modules over $\p[[T]]$: an element of a module over $\p[[T]]$ will be called "small" if its coefficients are small (a basis is assumed fixed).
\medskip
{\bf 1.4.3.} For $z=\sum_{i=0}^\infty \la_iT^i\in\p[[T]]$ we denote $z^{(k)}:=\sum_{i=0}^\infty \la_i^{q^k}T^i$.
$\tau$ acts on $\p\{T\}$ and $\p[[T]]$ by the formula $\tau(z)=z^{(1)}$, where $z\in \p\{T\}$ or $z\in \p[[T]]$.
\medskip
Now, we define (see (6.5) for the same constructions in a general setting):
$$M[[T]]:=M\otimes_{\p[T]}\p[[T]], \ \ M\{T\}:=M\otimes_{\p[T]}\p\{T\},\eqno{(1.4.3.1)}$$
$$M_{[[T]]}:=\Hom_{\p[T]}(M,\p[[T]]), \ \ M_{\{T\}}:=\Hom_{\p[T]}(M,\p\{T\})\eqno{(1.4.3.2)}$$
(see [G], 5.9.11.1 for (1.4.3.1)). We have: $\tau$ acts on $M[[T]]$, $M\{T\}$, $M_{[[T]]}$, $M_{\{T\}}$ by the standard formulas of the action of an operator on tensor products (see [G], 5.9.11.1) and Hom's.
\medskip
{\bf Definition 1.5.1.} $H^1(M)=M\{T\}^\tau=M[[T]]^\tau\cap M\{T\}$,
$$H_1(M)={M_{\{T\}}}^\tau=(\Hom_{\p[T]}(M,\p\{T\}))^\tau=\Hom_{\p[T,\tau]}(M,\p\{T\}).\eqno{(1.5.2)}$$
\medskip
Formula for $H^1(M)$ is [G], (5.9.11.2), formula for $H_1(M)$ is [G], (5.9.11.3), where $E$ of [G], (5.9.11.3) is $E(M)$, see (0.2). Formula (1.5.2) follows immediately from [G], (5.9.25). Both $H^1(M)$, $H_1(M)$ are free $\n F_q[T]$-modules.
\medskip
It is convenient to consider a more general object. Namely, let $(V, V_0, \tau)$ (abbreviations: $V$ or $(V,\tau)$) be a triple consisting of a free $\p[T]$-module $V_0$ of finite rank $r$,  a map $\tau: V_0 \to V_0$ which is a set isomorphism satisfying $\tau(gv_0)=g^{(1)}\tau(v_0)$, where $g\in \p[T]$, $v_0\in V_0$, and a fixed isomorphism $V_0\otimes_{\p[T]}\p[[T]]\to V$. Particularly, $V$ is a free $\p[[T]]$-module. We extend the action of $\tau$ on $V$ by the formula $\tau(v_0\otimes g)= \tau(v_0)\otimes g^{(1)}$, where $g\in \p[[T]]$, $v_0\in V_0$.
\medskip
Let $v_*=(v_1, \dots, v_r)^t$ be a basis of $V_0$. We define $Q$ analogously to (1.4.1), i.e. $$Qv_*=\tau(v_*).\eqno{(1.5.3)}$$

We define a subset $V_s$ ("s" means "small") of $V$ as a set of elements of the form $\sum_{i=0}^rf_iv_i$ where $f_i\in \p\{T\}$. It is a free $\p\{T\}$-module of rank $r$, it does not depend on a choice of $v_*$. Further, we define $V^\tau$, $V_s^\tau:=(V_s)^\tau=(V^\tau)_s$ as the sets of $\tau$-invariant elements. We have: $V^\tau$ is a $\n F_q[[T]]$-module and $V_s^\tau$ is a $\n F_q[T]$-module. We prove in (6.7.1) that $\dim_{\n F_q[[T]]}V^\tau=r$.
\medskip
{\bf Proposition 1.5.4.} In the above notations, $\dim_{\n F_q[T]}V_s^\tau\le r$.
\medskip
{\bf Proof.} This is [P], Lemma 3.3.7, with minor changes. Let us repeat the proof. We denote (like in [P]) the fraction field of $\p\{T\}$ by $\n L$.
\medskip
{\bf Lemma 1.5.4.1.} If $\mu_1, \dots, \mu_k\in V_s^\tau$ are linearly independent over $\n F_q[T]$ then they are linearly independent over $\p\{T\}$.
\medskip
{\bf Proof.} Let us assume the opposite: $\exists$ \ $f_i\in \p\{T\}$, not all 0, such that $\sum_{i=1}^k f_i\mu_i=0$, and let $\mu_1, \dots, \mu_k$ be a minimal set of elements having this property (i.e. all such sets of $k-1$ elements are linearly independent over $\p\{T\}$). We consider the same equality in $V_s\otimes_{\p\{T\}}\n L$. In this space we can assume that $f_1=1$ and other $f_i\in \n L$. Applying $\tau$ we get $\sum \tau(f_i)\mu_i=0$. Subtracting these equalities we get $\sum_{i=2}^k (f_i-\tau(f_i))\mu_i=0$. Since $\mu_1, \dots, \mu_k$ is a minimal linearly dependent over $\p\{T\}$ system, we get that $f_2, \dots, f_k\in \n L^\tau$. We have $\n L^\tau=\n F_q(T)$ ([P], Lemma 3.3.2), hence the coefficients of the sum $\mu_1+\sum_{i=2}^k f_i\mu_i=0$ are in $\n F_q(T)$. Multiplying by a common denominator we get a non-trivial linear dependence relation in $\n F_q[T]$ - a contradiction. $\square$
\medskip
The proposition follows immediately. Really, since $\p\{T\}$ is a principal ideal domain ([P], below (2.2.4.1)) and $V_s$ is a free $\p\{T\}$-module of rank $r$, we get that $V_s^\tau$ cannot contain more than $r$ elements linearly independent over $\n F_q[T]$. $\square$
\medskip
{\bf Corollary 1.5.5.} $h_1(M)\le r$, $h^1(M)\le r$. $\square$
\medskip
Let us return to the case $V=M[[T]]$ (i.e. $V_s=M\{T\}$). Using the above basis $f_*$ we can identify row vectors $Y\in M_{1\times r}(\p\{T\})= \p\{T\}^r$ with elements of $M\{T\}$: $Y \mapsto Y\cdot f_* \in M\{T\}$ where $Y\cdot f_*$ is the product of $1\times r$ and $r\times 1$ matrices. Also, we can identify column vectors $X\in \p\{T\}^r$ with elements of $\Hom_{\p[T]}(M,\p\{T\})$. Namely, let $\vf\in \Hom_{\p[T]}(M,\p\{T\})$. The column vector $X$ corresponding to $\vf$ is, by definition, $\left(\matrix \vf(f_1)\\ \dots\\ \vf(f_r) \endmatrix \right)$. Under these identifications, (1.4.1) gives us immediately
$$Y\in H^1(M)\iff Y^{(1)}Q=Y;\eqno{(1.6.1)}$$
$$X\in H_1(M)\iff QX=X^{(1)}.\eqno{(1.6.2)}$$

For a general case of arbitrary $(V,\tau)$ we have similar formulas. Let $Q$ be from (1.5.3), and $Y\in M_{1\times r}(\p\{T\})$ or $Y\in M_{1\times r}(\p[[T]])$ a row vector. The condition $Y\cdot v_*\in V^\tau$ is equivalent to
$$Y^{(1)}Q=Y, \hbox{ where } Y\in M_{1\times r}(\p[[T]])\eqno{(1.7.1)}$$
and the condition $Y\cdot v_*\in V_s^\tau$ is equivalent to
$$Y^{(1)}Q=Y, \hbox{ where } Y\in M_{1\times r}(\p\{T\}).\eqno{(1.7.2)}$$

{\bf 1.7.3.} The definition of the pairing $\pi$ from (0.2.1) is straightforward: if $y=\sum m_i\otimes Z_i\in H^1(M)$, where $m_i\in M$, $Z_i\in \p\{T\}$ and $\vf: M\to\p\{T\}$ belongs to $H_1(M)$ then $\pi(y\otimes\vf)=\sum \vf(m_i)Z_i$. In coordinates, $\pi$ is given by the following formula:
$$\pi(X,Y)=YX.\eqno{(1.7.4)}$$
Really, we have $(YX)^{(1)}=YQ^{-1}QX=YX$, i.e. $\pi(y\otimes\vf)\in \n F_q[T]$.
\medskip
{\bf Remark 1.7.5.} The proof of Theorem 0.2.2 given by Anderson is rather complicated. In reality, its part affirming $h^1(M)=r \iff h_1(M)=r$ has an easy proof. Let us give it.
\medskip
{\bf 1.7.5.1.} First, we consider an equation where $Z\in \p\{T\}$ is an unknown:
$$Z=(T-\th)Z^{(1)}.\eqno{(1.7.5.2)}$$
There exists a solution to (1.7.5.2) $Z=\Xi$ where $\Xi=\sum_{i=0}^\infty\xi_iT^i$ is from [G], Example 5.9.36, p. 172; it satisfies $\xi_0\ne0$, it is defined uniquely up to a factor from $\n F_q^*$ (also, $\Xi=1/\omega_1$ where $\omega_1$ is the Anderson-Thakur function defined in [AT], Section 2.5). We have: all roots to (1.7.5.2) have the form $U\cdot\Xi$ where $U\in \n F_q[T]$. Further, we have $\Xi^{-1} \in \p\{T\}$.
\medskip
Let $h^1(M)=r$. We consider a $r\times r$ matrix $\g Y\in GL_r(\p\{T\})$ formed by rows $Y_i$, $i=1,\dots,r$, corresponding to a $\n F_q[T]$-basis of $H^1(M)$, according to (1.6.1). Its lines are linearly independent over $\n F_q[T]$, hence, according Lemma 1.5.4.1, we get that they are linearly independent over $\p\{T\}$. This means that $\det \g Y\ne0$, $\det \g Y\in \p\{T\}$. Equation (1.6.1) applied to $\g Y$ becomes of the form
$$\g Y^{(1)}Q=\g Y.\eqno{(1.7.5.3)}$$
Taking the determinant of (1.7.5.3) and taking into consideration Lemma 1.4.2 we get $\det\g Y^{(1)}\cdot c \cdot (T-\th)^n=\det\g Y$. (1.7.5.1) shows that $\det \g Y\in\bar \n F_p\cdot \n F_q[T]\cdot\Xi^n$.
\medskip
Now we define $\g X:=\g Y^{-1}$. Since $\Xi^{-1} \in \p\{T\}$, we get: $\g X\in \n F_q(T)\cdot GL_r(\p\{T\})$. It satisfies $Q\g X=\g X^{(1)}$. After multiplication by a non-zero element of $\n F_q[T]$ we get $\g X\in GL_r(\p\{T\})$, it continue to satisfy $Q\g X=\g X^{(1)}$. This means that the columns of $\g X$ form $r$ column vectors satisfying (1.6.2), hence $h_1(M)=r$.
\medskip
These considerations are convertible, i.e. $h_1(M)=r\implies h^1(M)=r$.
\medskip
{\bf 1.8.} We shall use the theory of the dual t-motives. There are different notions of "dual", "duality" for t-motives, let us describe them. In the present paper we shall use only the notion of duality defined in [GL07]; this dual of $M$ is defined by $M'$.
\medskip
{\bf 1.8.1.} Duality in the meaning of [HJ], [P]. The main objects of [HJ] are A-motives ([HJ], Definition 3.1). They are a version of t-motives of the present paper. More exactly, an Anderson t-motive of the present paper is an A-motive of [HJ], it is called an effective A-motive ([HJ], Definition 3.1.c).
\medskip
For all A-motives $M$ of [HJ] there exists the dual A-motive $M^{\vee}$ in the meaning of [HJ]. The functor $M\mapsto M^{\vee}$ is an involutive auto-anti-equivalence of the category of A-motives.
\medskip
The relation between two notions of duality is the following. Let $M$ be an Anderson t-motive of the present paper having $N=0$; we denote the corresponding effective A-motive of [HJ] by $HJ(M)$. Let $\g C$ be the Carlitz module (see Remark 1.2.2).
\medskip
$M'$ (meaning of [GL07]) exists iff $(HJ(M))^{\vee}\otimes \g C$ is an effective A-motive of [HJ]; in this case $$HJ(M')=(HJ(M))^{\vee}\otimes \g C.\eqno{(1.8.1.1)}$$

Practically the same notion of duality is defined in [P], Definition 3.2.10; the same notation $P^{\vee}$ for the dual of $P$ is used.
\medskip
{\bf 1.8.2.} Let $M$ be an Anderson t-motive and $E(M)$ the corresponding t-module. Some authors call $E(M)$ the dual of $M$, because the functor $M \mapsto E(M)$ is anti-equivalent.
\medskip
{\bf 1.8.3.} Maurischat in [M], very beginning of Section 4 considers a t-module $E$ and an associated object $\g M$, he calls it a dual t-motive. $\g M$ is covariant with respect to $E$, hence it is contravariant with respect to a t-motive $M$ such that $E=E(M)$. This $\g M$ has nothing common with both above $M'$, $M^{\vee}$. It is a $\p\{\s\}$-module (where $\s=\tau^{-1}$) but not a $\p\{\tau\}$-module.
\medskip
Therefore, we see that only $M^{\vee}$ is related with $M'$ by (1.8.1.1), while $E(M)$ and $\g M$ are not related with $M'$.
\medskip
Recall that later the word "dual" will always mean $M'$ --- the dual in the meaning of [GL07]. Its disadvantage is that $M'$ is not defined for all $M$, but there are more advantages: formulas (1.9.1), lemma 1.10.2, commutativity with the lattice functor, close analogy with the duality of abelian varieties with MIQF ([GL09]) and of abelian varieties of CM-type ([GL07], Section 12, especially Remark 12.4.1).
\medskip
$M'$ is defined in [GL07], (1.8); for an explicit formula see [GL07], (1.10.1). Namely, let $M$ be a t-motive and its $Q$ from (1.4.1). If there exists a t-motive $M'$ such that its $Q(M')$ satisfies
$$Q(M')=(T-\th)(Q^t)^{-1}\eqno{(1.8.4)}$$

then $M'$ is called the dual of $M$. If it exists then it is well-defined, see [GL07].
\medskip
Not all Anderson t-motives have dual. We shall use only t-motives defined by the formula (1.10.1) with $N=0$, they all have dual, see Lemma 1.10.2. It is known that all pure t-motives having $N=0$ have a dual, except the trivial exception of t-motives having $r=n$ ([GL07], Theorem 10.3). Explicit formulas for duals of a large class of t-motives are given in [GL07], Section 11.
\medskip
The only application of the notion of duality that we need is the Proposition 1.9, it permits to reduce the calculation of $H_1(M)$ to the calculation of $H^1(M')$.
\medskip
{\bf Proposition 1.9.} If the dual t-motive $M'$ exists then there exists a canonical isomorphism
$$H^1(M)\to H_1(M')\eqno{(1.9.1)}$$
and hence $H_1(M)\to H^1(M')$, because $(M')'=M$.
\medskip
{\bf Proof.} We identify $H^1$, resp. $H_1$ with $Y$, $X$ as above. Hence, we get: $Y$ is a root to (1.6.1) for $M$  $\iff \ \Xi^{-1}Y^t$ is a root to (1.6.2) for $M'$. This gives a formula for the map (1.9.1) in coordinates: $Y\mapsto \Xi^{-1}Y^t$. It is easy to see that the map is canonical, i.e. it does not depend on a choice of a basis. $\square$
\medskip
We shall use this proposition in order to deal with $H^1(M')$ instead of $H_1(M)$.
\medskip
{\bf 1.10.} We shall find a counterexample to (0.3a) among t-motives defined by the equation (1.3) such that $l=2$, $\g A_2=I_n$. We denote $\g A_1$ by $A$, hence (1.3) has the form $$Te_*=(\theta I_n+N) e_*+A\tau e_* +\tau^2e_*.\eqno{(1.10.1)}$$ This t-motive is denoted by $M(A)=M(A,N)$ (a high-dimensional Drinfeld module of high-dimensional rank 2).\footnotemark \footnotetext{This is the terminology of the authors. Namely, an Anderson t-motive given by the formula $Te_*=\sum_{i=0}^k A_i\tau^ie_*$ is called a high-dimensional Drinfeld module of high-dimensional rank $k$ if $\det A_k\ne0$ (particularly, if $A_k=I_n$). The name "module" represents the same abuse of language as in Remark 1.2.2: really, it is a t-motive of rank $r=nk$.}
The rank of $M(A)$ is $2n$, hence, if $M(A)$ is uniformizable then its lattice (it is contained in $\p^n$) is isomorphic to $(\n F_q[\th])^{2n}$, i.e. they are the nearest analogs of abelian varieties of dimension $n$. A basis $f_*$ can be chosen as $(e_1, \dots, e_n, \tau(e_1), \dots, \tau(e_n) \ )^t$. The matrix $Q$ of $M(A)$ in this basis is $\left(\matrix 0&&I_n\\(T-\th)I_n-N&&-A \endmatrix \right)$ (entries are $n\times n$-blocks).
\medskip
Taking into consideration Proposition 1.9, we shall calculate $h^1((M(A))')$ instead of $h_1(M(A))$. So, we need
\medskip
{\bf Lemma 1.10.2.} If $N=0$ then $(M(A))'=M(A^t)$.
\medskip
{\bf Proof.} This is [GL07], formulas (7.1) --- (7.2), and Section 11. Since [GL07], Section 11 gives formulas for much more general situation, we give here an explicit proof. We have $Q((M(A))')=(T-\th)\left(\matrix 0&I_n\\(T-\th)I_n&-A \endmatrix \right)^{t-1}=\left(\matrix A^t&(T-\th)I_n \\ I_n &0 \endmatrix \right)$ (here $\g Q^{t-1}$ means $(\g Q^t)^{-1}$). This is $Q(M(-A^t))$ in the basis $(\tau(e_1),\dots, \tau(e_n), e_1, \dots, e_n)^t$. Since t-motives $M(A)$, $M(-A)$ are isomorphic (the isomorphism can be got by a scalar change of basis of $M(A)$ over $\p\{\tau\}$ with coefficient $\eta$ satisfying $\eta^{q-1}=-1$), we get the desired result. $\square$
\medskip
{\bf 1.10.3.} Hence, in order to prove that not always $h^1(M)=h_1(M)$, it is sufficient to find a matrix $A$ such that $h^1(M(A))\ne h^1(M(A^t))$.
\medskip
\medskip
\medskip
{\bf 2. Affine equations.}
\medskip
{\bf 2.0.1.} We consider $(V,\tau)$ as above. In order to find $\dim_{\n F_q[T]}V_s^\tau$ we must solve (1.7.1) and calculate the quantity of linearly independent small solutions (i.e. solutions to (1.7.2)). We have: $Q$ from (1.5.3) is fixed and $Y=(y_1,\dots,y_r)$ are unknowns. We denote $y_i=\sum _{j=0}^\infty y_{ij}T^j$. We have: (1.7.1) is a system of algebraic equations in $y_{ij}\in\p$, $i=1,\dots, r$, $j=0, \dots, \infty$. Conjecturally, for many $(V,\tau)$ we can transform (1.7.1) to the below system of equations \{E\} of type (2.1) in $y_{i_0,j}$ where $i_0\in [1,\dots, r]$ is fixed and $j=0, \dots, \infty$.
\medskip
{\bf 2.0.2.} Here "transform" means that if $y_{ij}$, $i=1,\dots, r$, $j=0, \dots, \infty$ is a solution to (1.7.1) then $y_{i_0,j}$, $j=0, \dots, \infty$ is a solution to \{E\}, and conversely, if $y_{i_0,j}$, $j=0, \dots, \infty$ is a solution to \{E\}  then there exists the only $y_{ij}$, $i=1,\dots, r$, $j=0, \dots, \infty$ which is a solution to (1.7.1).
\medskip
Examples of such transformations for $(V,\tau)$ coming from some $M$ (i.e. $V=M[[T]]$) are given in Sections 3, 6.1, 6.4. Therefore, let us give some general definitions and elementary results concerning such equations.
\medskip
Let $r,n\ge1$, $\vk_1,\dots,\vk_n\ge0$ be degrees, $a_\ga\in \p$  ($\ga=0,\dots,r$) and $b_{\be\ga}\in \p$ ($\be=1,\dots,n, \ \ \ga=0,\dots,\vk_\be$) coefficients, and $x_0,x_1,x_2,...\in \p$ unknowns. An $i$-th affine equation ($i=0,1,2,\dots$) has the form (here $x_j=0$ for $j<0$):

$$\sum_{\ga=0}^r a_\ga x_i^{q^\ga}+\sum_{\be=1}^n\sum_{\ga=0}^{\vk_\be} b_{\be\ga} x_{i-\be}^{q^\ga}=0.\eqno{(2.1)}$$

The set of terms $a_\ga x_i^{q^\ga}$ is called the head of the equation, the set of other terms is called the tail of the equation. Explicitly:
$$\hbox{The 0-th equation:} \ \ \ \ a_rx_0^{q^r}+a_{r-1}x_0^{q^{r-1}}+...+a_1x_0^q+a_0x_0=0\ \ \hbox{(head)}; \eqno{(2.1.0)}$$
$$\hbox{The first equation:} \ \ \ \ a_rx_1^{q^r}+a_{r-1}x_1^{q^{r-1}}+...+a_1x_1^q+a_0x_1+\ \ \hbox{(head)} \eqno{(2.1.1)}$$ $$ +b_{1,\vk_1}x_0^{q^{\vk_1}}+ b_{1,\vk_1-1}x_0^{q^{\vk_1-1}}+...+b_{11}x_0^q+b_{10}x_0=0\ \ \hbox{(tail)};$$
$$\hbox{The second equation:} \ \ \ \ a_rx_2^{q^r}+a_{r-1}x_2^{q^{r-1}}+...+a_1x_2^q+a_0x_2+\ \ \hbox{(head)} \eqno{(2.1.2)}$$ $$ +b_{1,\vk_1}x_1^{q^{\vk_1}}+ b_{1,\vk_1-1}x_1^{q^{\vk_1-1}}+...+b_{11}x_1^q+b_{10}x_1+\ \ \hbox{(tail)}$$ $$ +b_{2,\vk_2}x_0^{q^{\vk_2}}+ b_{2,\vk_2-1}x_0^{q^{\vk_2-1}}+...+b_{21}x_0^q+b_{20}x_0=0;$$
etc., the length of the tail\footnotemark \footnotetext{The length of the tail is the quantity of unknowns in the tail part of the equation, i.e. it is 0, respectively 1, 2 in (2.1.0), respectively (2.1.1), (2.1.2) etc.} stabilizes, which does not exceed $n$.

The equations of type (2.1) are called the equations of bounded tail, $n$ is its (maximal) length. We can consider the affine equations "with unbounded tail" as well ("u" means unbounded):

$$\sum_{\ga=0}^r a_\ga x_i^{q^\ga}+\sum_{\be=1}^\infty\sum_{\ga=0}^{\vk_\be} b_{\be\ga} x_{i-\be}^{q^\ga}=0.\eqno{(2.1u)}$$

Without loss of generality we shall consider equations satisfying $a_0\ne0$ (separability), $a_r\ne0$, $\forall \ \be \ \ \vk_\be<r$, $b_{\be,\vk_\be}\ne0$. Also, we can assume $a_r=1$. A substitution $x_i \mapsto \la x_i$ gives a change of coefficients $a_\ga \mapsto \la^{-q^\ga}a_\ga, \ \ b_{\be\ga} \mapsto \la^{-q^\ga}b_{\be\ga}$.
%Are there more isomorphisms of equations?.
\medskip
The systems (2.1), (2.1u) are solved consecutively: for $i=0$ the tail is 0, the set of $x_0$ --- solutions to (2.1.0) --- is a $\n F_q$-vector subspace of $\p$ of dimension $r$, denoted by $S_0$. Let $x_0, x_1,\dots,x_\al$ be a solution to (2.1), (2.1u) for $i=0,1,\dots, \al$. For $i=\al+1$ the equations (2.1), (2.1u) have the form

$$\sum_{\ga=0}^r a_\ga x_{\al+1}^{q^\ga}+W=0\eqno{(2.2.1)}$$
where $W=W(\al;x_{0},\dots,x_\al) \in\p$ is obtained by substitution of $x_{0},\dots,x_\al$ to the tail members. The set of solutions to (2.2.1) (for $x_0, x_1,\dots,x_\al$ fixed) is an affine space over $\n F_q$, with the base vector space $S_0$ (this explains the terminology).
\medskip
Let $x_0, x_1,\dots$ be a solution to (2.1), (2.1u). We associate it an element $\{x\}:=\sum_{i=0}^\infty x_iT^i\in \p[[T]]$ which (by abuse of language) will be also called a solution to (2.1), (2.1u).
\medskip
$\p[[T]]$ is a left $\p[T,\tau]$-module, the multiplication by $\tau$ is defined by (1.4.3). We can consider the multiplication by elements of $\p[T,\tau]$ as an action of an operator on $\p[[T]]$. From this point of view, we can consider (2.1) as an equation (here $P\in \p[T,\tau]$):
$$P(\{x\})=0\eqno{(2.2.2)}$$ where $\{x\}$ is as above and $P:=\sum_{\ga=0}^r a_\ga\tau^\ga + \sum_{\be=1}^n\sum_{\ga=0}^{\vk_\be} b_{\be\ga} \tau^\ga T^\be$.
\medskip
We denote the set of solutions to (2.1), (2.1u) in $\p[[T]]$ by $X=X(a_*,b_{**})$.
\medskip
We need a definition.
\medskip
{\bf Definition 2.2.3.} Let $\g W$ be a $\n F_q[[T]]$-submodule of $\p[[T]]$. $\g W$ is called a $T$-divisible submodule if for any $y\in \g W\cap T\cdot \p[[T]]$ we have $y/T\in \g W$.
\medskip
If we denote an element of $\g W$ by $y_i$ then we denote its coefficients by $y_{ij}$, i.e. $y_i=\sum_{j=0}^\infty y_{ij}T^j$. For any $\n F_q[[T]]$-submodule $\g W$ of $\p[[T]]$ we denote by $\g W(0)$ the set of free terms of elements of $\g W$. It is a $\n F_q$-vector subspace of $\p$.
\medskip
For any affine equation (2.1), (2.1u) we have: the set of its solutions $X$ is a free $T$-divisible $\n F_q[[T]]$-submodule of $\p[[T]]$, and $X(0)$ is $S_0$ in the above notations.
\medskip
{\bf Proposition 2.3.} Let $\g W$ be a $T$-divisible $\n F_q[[T]]$-submodule of $\p[[T]]$.
\medskip
(A) If $\g W$ is a finitely generated $\n F_q[[T]]$-module and the elements $y_1,\dots, y_k$ form a basis of $\g W$ over $\n F_q[[T]]$ then $\g W(0)$ is finitely-dimensional over $\n F_q$ and $y_{10},\dots,y_{k0}$ form a $\n F_q$-basis of $\g W(0)$.
\medskip
(B) Conversely, let $\g W(0)$ be finitely-dimensional over $\n F_q$ and $z_1,\dots, z_k$ form a $\n F_q$-basis of $\g W(0)$. Let $\forall \ i$ we choose elements $y_i\in \g W$ such that $z_i=y_{i0}$. Then $\g W$ is a finitely generated $\n F_q[[T]]$-module, and the elements $y_1,\dots, y_k$ form a basis of $\g W$ over $\n F_q[[T]]$.
\medskip
{\bf Proof.} (A) The fact that $y_{10},\dots,y_{k0}$ generate $\g W(0)$ over $\n F_q$ is obvious. Let us assume that they are dependent over $\n F_q$: $\exists \ c_i \in \n F_q$ such that $\sum_{i=1}^n c_i y_{i0}=0$. In this case $\sum_{i=1}^n c_i y_{i}\in \g W\cap T\cdot \p[[T]] $. Since $\g W$ is $T$-divisible, we have $(\sum_{i=1}^n c_i y_{i})/T\in \g W$ and hence (because $y_1,\dots, y_k$ generate $\g W$) \ $\exists \ g_i\in \n F_q[[T]]$ such that
$$(\sum_{i=1}^n c_i y_{i})/T=\sum_{i=1}^n g_i y_{i}.$$
This means that $\sum_{i=1}^n (c_i-Tg_i) y_{i}=0$ is a non-trivial dependence relation on $y_i$ --- a contradiction.\footnotemark \footnotetext{[M], proof of Proposition 3.4 (2) is similar to the proof of the present proposition.}
\medskip
(B). Let us assume that $y_1,\dots, y_k$ are dependent, i.e. \ $\exists \ \ g_i\in \n F_q[[T]]$ such that $\sum_{i=1}^n g_i y_{i}=0$. Dividing $g_i$ by a power of $T$ we can assume that some $g_{i0}$ are non-0. We get a non-trivial relation $\sum_{i=1}^n g_{i0} z_{i}=0$ --- a contradiction to linear independence of $z_i$.
\medskip
Let us prove that $y_1,\dots, y_k$ generate $\g W$ over $\n F_q[[T]]$. Let $w_0\in \g W$. We have: $\forall \ i \ \exists \ c_{i0}\in \n F_q$ such that $w_{00}= \sum_{i=1}^n c_{i0} z_{i}$ and hence $w_0-\sum_{i=1}^n c_{i0} y_{i}\in T\cdot \n F_q[[T]] \cap \g W$. We denote $w_1:=(w_0-\sum_{i=1}^n c_{i0} y_{i})/T$. Since $\g W$ is $T$-divisible, we have $w_1\in \g W$. Now we find $c_{i1}\in \n F_q$ such that $w_{10}= \sum_{i=1}^n c_{i1} z_{i}$ and denote $w_2:=(w_1-\sum_{i=1}^n c_{i1} y_{i})/T$. Further, we find find $c_{i2}\in \n F_q$ such that $w_{20}= \sum_{i=1}^n c_{i2} z_{i}$ etc. Continuing this process we get $c_i:=\sum_{j=0}^\infty c_{ij}T^j$ and $w_0= \sum_{i=1}^n c_i y_i$. $\square$
\medskip
{\bf Corollary 2.3.1.} $X$ (the set of solutions to (2.1), (2.1u)) is a free $\n F_q[[T]]$-module of rank $r$.
\medskip
The converse is also true:
\medskip
{\bf Proposition 2.3.2.} Let $\g W$ be a finitely generated $T$-divisible $\n F_q[[T]]$-submodule of $\p[[T]]$, and let $\dim \g W(0)=r_0$. There exists the only affine equation (2.1u) satisfying $r=r_0$, $a_r=1$, $a_0\ne0$, $\forall \ \be$ we have $\vk_\be< r$ and such that $\g W$ is the set of its roots.
\medskip
{\bf Proof.} We denote $\g P_0(x):= \sum_{\ga=0}^r a_\ga x^{q^\ga}$ (the head), and for all $\be>0$ we denote $\g P_\be(x):= \sum_{\ga=0}^{\vk_\be} b_{\be\ga} x^{q^\ga}$ (the $\be$-th polynomial of the tail).
\medskip
We have: $\g P_0(x)=\prod_{u\in \g W(0)}(x-u)$; it is a polynomial of the form (2.1.0) having $r=r_0$, $a_r=1$, $a_0\ne0$. Clearly it is the only possibility for $\g P_0$.
\medskip
Now we choose $y_1, \dots, y_r$ --- a basis of $\g W$ over $\n F_q[[T]]$. (2.1.1) means that $\forall \ i =1,\dots, r$ we must have $$\g P_1(y_{i0})=-\g P_0(y_{i1}).\eqno{(2.3.2.1)}$$ Let $C$ be a $r\times r$-matrix whose $(\al,\be)$-th entry is ${y_{\al0}}^{q^{\be-1}}$ (the transposed of the Moore matrix, see [G], Definition 1.3.2). We have: $|C|\ne0$ ([G], Lemma 1.3.3), hence for all numbers $\mu_1\dots,\mu_r$ there exists the only numbers $\ga_0\dots,\ga_{r-1}$ such that $\forall \ i =1,\dots, r$ we have
$$\sum_{j=0}^{r-1} \ga_j y_{i0}^{q^{j}}=\mu_i.\eqno{(2.3.2.2)}$$

Choosing $\mu_i=-\g P_0(y_{i1})$ we get that there exists the only polynomial $\g P_1(x)$ of the form $\g P_1(x)=\sum_{j=0}^{r-1} b_{1j}x^{q^j}$ such that (2.3.2.1) is satisfied.
\medskip
Further, (2.1.2) means that $\forall \ i =1,\dots, r$ we must have $$\g P_2(y_{i0})=-(\ \g P_0(y_{i2})+\g P_1(y_{i1})\ ).\eqno{(2.3.2.3)}$$

By the same argument we see that this $\g P_2(x)$ of the form $\g P_2(x)=\sum_{j=0}^{r-1} b_{2j}x^{q^j}$ really exists and is unique. Continuing this process we get the desired. $\square$
\medskip
This proposition permits to classify finitely generated $T$-divisible $\n F_q[[T]]$-submodules of $\p[[T]]$ according the length of the tail of the corresponding affine equation. So, a finitely generated $T$-divisible $\n F_q[[T]]$-submodule of $\p[[T]]$ can be of finite or unbounded tail; if it is of finite tail then $n$ \ --- \ the length of the tail \ --- \ is well-defined.
\medskip
Let $\g W$ be a finitely generated $T$-divisible $\n F_q[[T]]$-submodule of $\p[[T]]$. The elements of $\g W$ belonging to $\p\{T\}$ are called small elements. They form a $\n F_q[T]$-module. It is denoted by $\g W_s$.
\medskip
It is easy to prove that the $\n F_q[T]$-rank of $\g W_s$ is $\le$ the $\n F_q[[T]]$ -rank of $\g W$. First, we give a simple proof of this fact under a supposition that $\g W_s$ is finitely generated over $\n F_q[T]$:
\medskip
{\bf Proposition 2.3.3.} Let $y_1,\dots, y_k$ be a basis of $\g W_s$ over $\n F_q[T]$. Then $y_{1},\dots, y_{k}$ are linearly independent over $\n F_q[[T]]$.
\medskip
{\bf Proof.} This is practically the proof of (A) of 2.3, or, the same, the proof of [M], Proposition 3.4 (2). Let us repeat it. First, we prove that $y_{10},\dots, y_{k0}$ are linearly independent over $\n F_q$. We assume that $y_{10},\dots, y_{k0}$ are linearly dependent over $\n F_q$, namely, $\exists \ c_i \in \n F_q$ such that $\sum_{i=1}^k c_i y_{i0}=0$. We have $(\sum_{i=1}^k c_i y_{i})/T \in \g W_s$, hence $\exists \ g_i \in \n F_q[T]$ such that $(\sum_{i=1}^k c_i y_{i})/T=\sum_{i=1}^k g_i y_{i}$. This means that $\sum_{i=1}^k (c_i-Tg_i) y_{i}=0$ --- this is a non-trivial relation of dependency of $y_1,\dots, y_k$ over $\n F_q[T]$ --- a contradiction.
\medskip
Now we repeat the proof of (B) of 2.3. Namely, let us assume that $y_1,\dots, y_k$ are dependent over $\n F_q[[T]]$, i.e. \ $\exists \ \ g_i\in \n F_q[[T]]$ such that $\sum_{i=1}^n g_i y_{i}=0$. Dividing $g_i$ by a power of $T$ we can assume that some $g_{i0}$ are non-0. We get a non-trivial relation $\sum_{i=1}^n g_{i0} y_{i0}=0$ --- a contradiction to linear independence of $y_{i0}$ over $\n F_q$. $\square$
\medskip
Now we give a proof of this fact without supposition that that $\g W_s$ is finitely generated over $\n F_q[T]$:
\medskip
{\bf Proposition 2.3.4.} Let $f_0,\dots, f_n\in\g W_s$ be linearly independent over $\n F_q[T]$. Then they are linearly independent over $\n F_q[[T]]$.
\medskip
{\bf Proof.} We need a lemma:
\medskip
{\bf Lemma.} Let $f_1,\dots, f_n\in\p\{T\}$, $g_1,\dots, g_n\in\n F_q[[T]]$ such that $f_0:=\sum_{i=0}^ng_if_i\in\p\{T\}$. If $f_1,\dots, f_n$ are linearly independent over $\n F_q[T]$ then $g_1,\dots, g_n\in\n F_q(T)$.
\medskip
{\bf Proof.} Induction by $n$. Let $n=1$. Multiplying $f_1$ by a scalar we can assume that all coefficients of $f_1$ have $\ord\ge0$, and some of them have $\ord=0$. We denote the reduction $O_{\p}\to \bar \n F_q$ by bar. If $f_i\in O_{\p}\{T\}$ then $\bar f_i\in \bar \n F_q[T]$. We have $g_1f_1=f_0$; reducing this equation we have $g_1\bar f_1=\bar f_0$ where $\bar f_1\ne0$. Hence, $g_1=\bar f_0/\bar f_1\in \bar \n F_q(T)$. Since $\bar \n F_q(T)\cap \n F_q[[T]]\subset \n F_q(T)$, the lemma for $n=1$ is proved.
\medskip
Now we assume that the lemma is proved for $n-1$; let us prove it for $n$.
Again multiplying all $f_i$ by a scalar we can assume that all coefficients of all $f_i$ have $\ord\ge0$ and $\exists \ i=1,\dots,n$ such that $\bar f_i\ne 0$. Reducing $$\sum_{i=1}^ng_if_i=f_0\eqno{(2.3.4.1)}$$ we get
$$\sum_{i=1}^ng_i\bar f_i=\bar f_0.\eqno{(2.3.4.2)}$$
$\exists \ m$ such that $\forall \ i \ \bar f_i\in \n F_{q^m}[T]$. Let $\al_1, \dots, \al_m$ be a basis of $\n F_{q^m}$ over $\n F_{q}$, and let
$$\bar f_i=\sum_{j=1}^m h_{ij}\al_j,$$
where $h_{ij}\in \n F_{q}[T]$.
We can assume $\bar f_n\ne0$ and $h_{n1}\ne0$. The $\al_1$-th coordinate of (2.3.4.2) becomes
$$\sum_{i=1}^ng_ih_{i1}=h_{01}.\eqno{(2.3.4.3)}$$

Multiplying (2.3.4.1) by $h_{n1}$ and (2.3.4.3) by $f_n$ we get
$$\sum_{i=1}^ng_if_ih_{n1}=f_0h_{n1},\eqno{(2.3.4.4)}$$
$$\sum_{i=1}^ng_ih_{i1}f_n=h_{i0} f_n.\eqno{(2.3.4.5)}$$

Subtracting (2.3.4.5) from (2.3.4.4) we get
$$\sum_{i=1}^{n-1}g_i(f_ih_{n1}-h_{i1}f_n)=f_0h_{n1}-h_{01}f_n.\eqno{(2.3.4.6)}$$

All $f_ih_{n1}-h_{i1}f_n$ are small. Further, if $f_1, \dots, f_n$ are linearly independent over $\n F_q[T]$ then $f_ih_{n1}-h_{i1}f_n$ are linearly independent over $\n F_q[T]$. By induction supposition, we get $g_1, \dots, g_{n-1}\in \n F_q(T)$. Since $\bar f_n\ne0$, (2.3.4.3) implies $g_{n}\in \n F_q(T)$. $\square$
\medskip
Now we deduce the proposition from the lemma. Let us assume the opposite: $\exists \ f_0,\dots, f_n\in\p\{T\}$ linearly independent over $\n F_q[T]$ and linearly dependent over $\n F_q[[T]]$. Therefore, $\exists \ g_0,\dots, g_n\in\n F_q[[T]]$ such that $\sum_{i=0}^ng_if_i=0$. Dividing if necessary all $g_i$ by a power of $T$ and renumbering them we can assume that the free term of $g_0$ is $\ne0$. Dividing by $g_0$ and changing notations we can assume $$\sum_{i=1}^ng_if_i=f_0,$$ $g_i\in\n F_q[[T]]$, $f_0,\dots, f_n, f\in\p\{T\}$ and $f_1,\dots, f_n$ are linearly independent over $\n F_q[T]$.
\medskip
The lemma implies that $g_1,\dots, g_n\in\n F_q(T)$. Multiplying by their common denominator we get a non-trivial dependence relation between $f_0,\dots, f_n$ over $\n F_q[T]$ --- a contradiction. $\square$
\medskip
Let us consider relations between $(V,\tau)$ and their $V^\tau$, $V_s$, $V_s^\tau$ of Section 1, and $T$-divisible $\n F_q[[T]]$-submodules of $\p[[T]]$ of the present section (see 2.0.1, 2.0.2). For $i=1,\dots, r$ let $\g p_i: \p[[T]]^r\to \p[[T]]$ be the projection to the $i$-th coordinate. We cannot guarantee that for all $(V,\tau)$ and $i$ we have: $\g p_i(V^\tau)$ is a $T$-divisible $\n F_q[[T]]$-submodule of $\p[[T]]$, we can only state a conjecture:
\medskip
{\bf Conjecture 2.3.5.} For all uniformizable $M$, for all $i\in \{1,\dots, r\}$, for $V=M[[T]]$ we have: $\g p_i(V^\tau)$ is a $T$-divisible $\n F_q[[T]]$-submodule of $\p[[T]]$.
\medskip
For $M$ defined by (1.10.1), $N=0$, $n=2$ this conjecture follows from (3.8) -- (3.10), for $n=3$ this is 6.4. For some $M$ having $n=2$, $r=5$ this is (6.2.3).
\medskip
Moreover, formulas (3.8) -- (3.10) and (6.2.3) suggest a stronger form of 2.3.5:
\medskip
{\bf Conjecture 2.3.6.} For the same $M$, $i$, $V$ we have: $\g p_i(V^\tau)$ is a $T$-divisible $\n F_q[[T]]$-submodule of $\p[[T]]$ having a finite tail, and $r$, $n$ of the tail are $r$, $n$ of $M$.
\medskip
{\bf Remark.} Coincidence of $r$ is obvious, see proof of Proposition 6.7.1, while coincidence of $n$ is mysterious.
\medskip
The below propositions of Section 2 are not necessary for the proof of Theorems 4.1, 4.7, 5.1. They are given for completeness and for future applications. For example, they were already used in [EGL].
\medskip
{\bf Definition 2.4.} Let $\{x\}=(x_0,x_1,\dots)$ be a solution to (2.1). It is called simple (or of simple type) if for all $i_0$ we have: ord's of all tail members of the equation ((2.1), $i=i_0$) for this $\{x\}$ (i.e. obtained while we substitute $x_0,x_1,...,x_{i_0-1}$) are different. An equation (2.1) is called simple if all its solutions are simple.
\medskip
Particularly, if the tail contains one term then the equation is trivially simple.
\medskip
For simple equations we can easily find $\ord x_i$. Really, we find all possible $\ord x_0$ treating the Newton polygon of the head of (2.1) for $i=0$. To pass from $i$ to $i+1$, we get that $\ord W$ (where $W$ is from (2.2.1)) is the minimum of the ord's of the tail terms. Again using the Newton polygon of the head of (2.1) and $\ord W$, we get ord's of all possible $x_{i+1}$.
\medskip
We see that the simplicity of (2.1) depends only on ord's of $a_\ga$, $b_{\be\ga}$. They belong to $\n Q\cup \infty$; for any $i$ the condition of non-simplicity imposes linear relations on ord's of $a_*$, $b_{**}$. We conjecture that for given $r$, $n$ there are only finitely many such relations, i.e. "almost all" equations are simple.
\medskip
Let $x_0$ be a fixed solution to (2.1) for $i=0$.
\medskip
{\bf Definition 2.5.} A solution $\{x\}=(x_0, x_1, x_2, \dots)$ to (2.1) is called a minimal chain generated by $x_0$ if it satisfies the following condition: $\forall \ i_0>0$ we have: $x_{i_0}$ is a solution to ((2.1), $i=i_0$) corresponding to the leftmost segment of the Newton polygon of ((2.1), $i=i_0$) for $x_0, x_1, x_2, \dots, x_{i_0-1}$ considered as parameters of ((2.1), $i=i_0$), i.e. $\ord x_{i_0}$ has the maximal possible value amongst ord's of solutions to ((2.1), $i=i_0$) for fixed $x_0, x_1, x_2, \dots, x_{i_0-1}$.
\medskip
A minimal chain generated by $x_0$ can be either simple or not.
If $x_0$ is fixed then a minimal chain generated by $x_0$ is not unique even if (2.1) is simple. But if there exists a simple minimal chain generated by $x_0$ (where (2.1) can be simple or not) then all minimal chains generated by $x_0$ are simple, and the sequence $\ord x_1, \ord x_2, \dots$ is uniquely defined by $\ord x_0$.

\medskip
{\bf Proposition 2.6.} Let $x_0\in S_0$ be such that its minimal chain $\{x\}=x_0, x_1, x_2, \dots$ is simple. Let $\{y\}=y_0, y_1, y_2, \dots$ be another simple solution to (2.1). Then $\ord y_0 \le \ord x_0$ (resp. $\ord y_0 < \ord x_0$) implies: $\forall \ i$ we have: $\ord y_i \le \ord x_i$ (resp. $\ord y_i < \ord x_i$).
\medskip
{\bf Proof.} Immediate, by induction. We consider the case $\ord y_0 \le \ord x_0$ (for the case $\ord y_0 < \ord x_0$ the proof is the same). Let this proposition be true for some $i$. We consider the equation (2.2.1) for $i+1$. We denote by $W_x$, resp. $W_y$ the term $W$ in (2.2.1) for the sets $x_0, x_1, x_2, \dots,x_i$, resp. $y_0, y_1, y_2, \dots,y_i$. Simplicity of $\{x\}$, $\{y\}$ implies that $\ord W_x, \ \ord W_y = $ minimum of the ord's of the corresponding terms of the tail. Hence, because $\ord y_j \le \ord x_j$ for $j=1,\dots,i$, we have $\ord W_y\le \ord W_x$, i.e. the leftmost vertex of the Newton polygon for (2.2.1) for $y_0, y_1, y_2, \dots,y_i$ is below or equal to the leftmost point of the Newton polygon for (2.2.1) for $x_0, x_1, x_2, \dots,x_i$. These two Newton polygons are the convex hulls of these points having $x$-coordinate 0 (the leftmost vertices), and other points corresponding to the head of (2.1) which are the same for $\{x\}$, $\{y\}$. This means that the inclination of the leftmost segment of the Newton polygon ( $=-\ord$ of the root) for (2.2.1) for $y_0, y_1, y_2, \dots,y_i$ (denoted by $\inc_{i+1}(y)$) is $\ge$ of the inclination of the leftmost segment of the Newton polygon for (2.2.1) for $x_0, x_1, x_2, \dots,x_i$ (denoted by $\inc_{i+1}(x)$).
\medskip
Inclinations of other sides of the Newton polygon for (2.2.1) for $y_0, y_1, y_2, \dots,y_i$ are $\ge \inc_{i+1}(y)$, hence $\ge \inc_{i+1}(x)$. This means that $\ord y_{i+1} \le \ord x_{i+1}$. $\square$
\medskip
{\bf Proposition 2.7.} Let $\{x\}$ be as above, and let $\{y\}=y_0, y_1, y_2, \dots$ be a minimal chain of $y_0$ (not necessarily simple). Then $\ord y_0 \ge \ord x_0$ (resp. $\ord y_0 > \ord x_0$) implies: $\forall \ i$ we have: $\ord y_i \ge \ord x_i$ (resp. $\ord y_i > \ord x_i$).
\medskip
{\bf Proof.} By induction, similar to the proof of Proposition 2.6. $\square$
\medskip
{\bf Corollary 2.8.} Let (2.1) be simple. If for minimal chains for all $x_0$ we have $lim_{i\to \infty}\ord x_i \ne +\infty$ then the dimension of (2.1) is 0.
\medskip
{\bf Conjecture 2.9.} Let (2.1) be simple, and let $x_{10}, x_{20}, \dots,x_{r0}$ be a $\n F_q$-basis of $S_0$. The dimension of (2.1) is equal to the quantity of $\al$ such that the minimal chain of $x_{\al0}$ is a small solution.
\medskip
{\bf Idea of the proof.} We can assume that $\ord x_{10} \ge \ord x_{20} \ge \dots \ge \ord x_{r0}$. Let $k$ be maximal number such that a minimal chain of $x_{k0}$ (denoted by $\{x_k\}$) is a small solution. According Proposition 2.6, we have that a minimal chain of $x_{\al0}$ is a small solution iff $\al\le k$. We must prove that minimal chains of $x_{10}, x_{20}, \dots,x_{k0}$ form a $\n F_q[T]$-basis of the set of small solutions. They are $\n F_q[[T]]$-, and hence $\n F_q[T]$-linearly independent. Let $\{y\}=(y_0, y_1,\dots)$ be a small solution. We apply Proposition 2.6 for $\{x\}=\{x_{k+1}\}$; it gives us that $y_0$ is a $\n F_q$-linear combination of $x_{10}, x_{20}, \dots,x_{k0}$, i.e. $y_0=\sum_{i=1}^k c_{i0}x_{i0}$ where $c_{i0}\in \n F_q$ are coefficients. We can consider $\{y\}(1):=(\{y\}-\sum_{i=1}^k c_{i0}\{x_i\})/T$ which is also a small solution. Applying the same operation to $\{y\}(1)$ we get $\{y\}(2)$ etc. As a result, we get that $\{y\}$ is a $\n F_q[[T]]$-linear combination of $x_{10}, x_{20}, \dots,x_{k0}$. We need to show that $\{y\}$ is a $\n F_q[T]$-linear combination of $x_{10}, x_{20}, \dots,x_{k0}$. This is an exercise for a student; we need this fact only for a case $k=1$, and $\forall \ i$ $\ord x_{1,i+1}>\ord x_{1i}$ where it is obvious (see proof of Lemma 4.6).
\medskip
Some cases of this conjecture were used in [EGL]; proofs for these cases were given by explicit calculations.
\medskip
There is a result for equations whose tail consists of one term:
\medskip
{\bf 2.11.} Let the only tail term be $b_{1k}x^{q^k}_{i-1}$ for some fixed $k$, and let $a_r=1$.
\medskip
We denote $\al_i:=\ord a_i$, $\be:=\ord b_{1k}$.
\medskip
{\bf Proposition 2.12.} Let 2.11 hold, and let $q^j$ be $x$-coordinate of the right end of the leftmost segment of the Newton polygon of the head of (2.1). Then
$$\frac{\al_0 -\al_j}{q^j-1}\le\frac{\al_0 -\be}{q^k-1}\eqno{(2.12.1)}$$ $$\iff \hbox{ the dimension of (2.1) is 0.}$$

{\bf Proof.} We denote $y_0:=\ord x_{10}$ where $x_{10}$ be a root to (2.1), $i=0$ corresponding to the leftmost segment of the Newton polygon of the head of (2.1). We have $y_0=\frac{\al_0 -\al_j}{q^j-1}$. Let us consider the equation (2.1), $i=1$ for this value of $x_{10}$. Ord of its free term is $\be+q^ky_0$. The negation of condition (2.12.1) is equivalent to $\be+q^ky_0-\al_0 > y_0$. Hence, if (2.12.1) does not hold then the leftmost segment of the Newton polygon of (2.1), $i=1$ is the segment $(0,\be+q^ky_0); (1,\al_0)$. Let $x_{11}$ be the root to to (2.1), $i=1$ corresponding to this segment. We denote $y_1:=\ord x_{11}=\be+q^ky_0-\al_0$. We have $y_1>y_0$. Hence, for $i=2$ the leftmost segment of the Newton polygon of (2.1), $i=2$ is the segment $(0,\be+q^ky_1); (1,\al_0)$. Continuing the process of finding the minimal chain corresponding to $x_{10}$ we get a solution $\sum_{j=0}^\infty x_{1j}T^j$. We denote $y_\ga:=\ord x_{1\ga}$, they satisfy a recurrent relation $y_{\ga+1}=\be+q^ky_\ga-\al_0$. A formula for $y_\ga$ is
$$y_\ga=\frac{\al_0 -\be}{q^k-1}+(\frac{\al_0 -\al_j}{q^j-1}-\frac{\al_0 -\be}{q^k-1})q^{k\ga},$$
it is proved immediately by induction. We get that if (2.12.1) does not hold then the dimension of (2.1) is $>0$.

Let us assume that (2.12.1) holds. In this case $\be+q^ky_0-\al_0 \le y_0$, hence $y_1\le y_0$. We get by induction that $\forall \ \ga $ we have $y_\ga\le y_0$, hence the minimal chain generated by $x_{10}$ is not small. Proposition 2.6 implies that the dimension of (2.1) is 0. $\square$
\medskip
{\bf 3. Affine equation corresponding to a t-motive.}
\medskip
Let $n=2$. We consider t-motives given by (1.10.1), where either $N=0$ or $N=N_0:=\left(\matrix 0&1\\ 0&0 \endmatrix \right)$, i.e. $N=\ve N_0$ where $\ve=0$ or 1. Let $A=\left(\matrix a_{11}& a_{12} \\  a_{21}& a_{22} \endmatrix \right)$.
\medskip
To find $h^1(M(A))$ we use formula (1.6.1). More exactly, after we eliminate 3 unknowns in (1.6.1), we get an affine equation and find its dimension. Namely, let $Y=(y_{11}, y_{12}, y_{21}, y_{22})$ be from (1.6.1) for $M(A)$. We denote $Y=(y_1, y_2)$ as a block matrix where $y_1=(y_{11}, y_{12})$, $y_2=( y_{21}, y_{22})$. Then (1.6.1) written in a 2-block form is\footnotemark \footnotetext{We give here all steps of these elementary calculations in order to simplify verification.}

$$(y_1^{(1)},y_2^{(1)})\left(\matrix 0&I_2\\(T-\th)I_2-\ve N_0&-A \endmatrix \right)=(y_1, y_2), \hbox{ i.e.}\eqno{(3.0a)}$$
$$y_2^{(1)}((T-\th)I_2-\ve N_0)=y_1, \ \ \ y_1^{(1)}-y_2^{(1)}A=y_2. \ \ \ \eqno{(3.0)}$$ Hence,
$$y_2=y_2^{(2)}((T-\th^q)I_2-\ve N_0)-y_2^{(1)}A.\eqno{(3.1)}$$
Substituting $y_2=(y_{21},y_{22})$ to (3.1) we get
$$y_{21}=y_{21}^{(2)}(T-\th^q)-y_{21}^{(1)}a_{11}-y_{22}^{(1)}a_{21};\eqno{(3.2)}$$
$$y_{22}=-\ve y_{21}^{(2)}+y_{22}^{(2)}(T-\th^q)-y_{21}^{(1)}a_{12}-y_{22}^{(1)}a_{22}.\eqno{(3.3)}$$ Now we eliminate $y_{22}$ from (3.2), (3.3).
Assuming $a_{21}\ne0$ we get from (3.2):

$$y_{22}^{(1)}=-\frac{1}{a_{21}}y_{21}+\frac{T-\th^q}{a_{21}}y_{21}^{(2)}-\frac{a_{11}}{a_{21}}y_{21}^{(1)}\eqno{(3.4)}$$ and hence
$$y_{22}^{(2)}=-\frac{1}{a_{21}^q}y_{21}^{(1)}+\frac{T-\th^{q^2}}{a_{21}^q}y_{21}^{(3)}-\frac{a_{11}^q}{a_{21}^q}y_{21}^{(2)};\eqno{(3.5)}$$
$$y_{22}^{(3)}=-\frac{1}{a_{21}^{q^2}}y_{21}^{(2)}+\frac{T-\th^{q^3}}{a_{21}^{q^2}}y_{21}^{(4)}-\frac{a_{11}^{q^2}}{a_{21}^{q^2}}y_{21}^{(3)}.\eqno{(3.6)}$$
From (3.3) we get $$y_{22}^{(1)}=-\ve y_{21}^{(3)}+y_{22}^{(3)}(T-\th^{q^2})-y_{21}^{(2)}a_{12}^q-y_{22}^{(2)}a_{22}^q.\eqno{(3.7)}$$ Substituting (3.4) - (3.6) to (3.7) we get
$$\frac{(T-\th^{q^3})(T-\th^{q^2})}{a_{21}^{q^2}}y_{21}^{(4)}+[(-\frac{a_{11}^{q^2}}{a_{21}^{q^2}}-\frac{a_{22}^q}{a_{21}^q})(T-\th^{q^2}) -\ve] y_{21}^{(3)}+$$ $$[-\frac{T-\th^q}{a_{21}}-\frac{T-\th^{q^2}}{a_{21}^{q^2}}+\frac{a_{11}^qa_{22}^q}{a_{21}^q}-a_{12}^q] y_{21}^{(2)} + (\frac{a_{11}}{a_{21}}+\frac{a_{22}^q}{a_{21}^q})y_{21}^{(1)}+\frac1{a_{21}}y_{21}=0.\eqno{(3.8)}$$ We denote $y_{21}$ by $x=\sum_{i=0}^\infty x_iT^i$ where $x_i\in\p$. Substituting this formula to (3.8) we get an affine equation of type (2.1) having $r=4$, $n=2$, $\vk_1=\vk_2=4$, and $a_\ga$, $b_{\be\ga}$ are the following:
\medskip
$$a_4=\frac{\th^{q^3+q^2}}{a_{21}^{q^2}}; \ \ \ a_3=\frac{a_{11}^{q^2}\th^{q^2}}{a_{21}^{q^2}}+\frac{a_{22}^q\th^{q^2}}{a_{21}^q}-\ve; \ \ \ a_2=\frac{\th^q}{a_{21}}+\frac{\th^{q^2}}{a_{21}^{q^2}}+\frac{a_{11}^qa_{22}^q}{a_{21}^q}-a_{12}^q; $$ $$ a_1=\frac{a_{11}}{a_{21}}+\frac{a_{22}^q}{a_{21}^q}; \ \ \ a_0=\frac1{a_{21}}; \eqno{(3.9)}$$ $$b_{14}=-\frac{\th^{q^3}+\th^{q^2}}{a_{21}^{q^2}}; \ \ \ b_{13}=-\frac{a_{11}^{q^2}}{a_{21}^{q^2}}-\frac{a_{22}^q}{a_{21}^q}; \ \ \ b_{12}=-\frac{1}{a_{21}}-\frac{1}{a_{21}^{q^2}}; \ \ \ b_{24}=\frac{1}{a_{21}^{q^2}}. $$

%\newpage
Hence, the equations \{E\} of type (2.1) (see (2.0.1)) for $M(A)$ have the form (we do not want to get $\forall \ \be \ \ \vk_\be<r$)
$$\frac{\th^{q^3+q^2}}{a_{21}^{q^2}}x_i^{q^4}+[(\frac{a_{11}^{q^2}}{a_{21}^{q^2}}+\frac{a_{22}^q}{a_{21}^q})\th^{q^2}-\ve]x_i^{q^3}+ (\frac{\th^q}{a_{21}}+\frac{\th^{q^2}}{a_{21}^{q^2}}+\frac{a_{11}^qa_{22}^q}{a_{21}^q}-a_{12}^q)x_i^{q^2}+ (\frac{a_{11}}{a_{21}}+\frac{a_{22}^q}{a_{21}^q})x_i^{q}+\frac1{a_{21}}x_i$$
$$-\frac{\th^{q^3}+\th^{q^2}}{a_{21}^{q^2}}x_{i-1}^{q^4}-(\frac{a_{11}^{q^2}}{a_{21}^{q^2}}+\frac{a_{22}^q}{a_{21}^q})x_{i-1}^{q^3}- (\frac{1}{a_{21}}+\frac{1}{a_{21}^{q^2}})x_{i-1}^{q^2}+\frac{1}{a_{21}^{q^2}}x_{i-2}^{q^4}=0.\eqno{(3.10)}$$
\medskip
{\bf Remark.} There exists another form to write (3.10):
$$[\th^{q^2}\tau^2+a_{22}^q\tau+1-\tau^2T][\frac1{a_{21}}(\th^{q}\tau^2+a_{11}\tau+1-\tau^2T)](x)=a_{12}^q\tau^2(x).\eqno{(3.11)}$$
which is much more "agreeable" than the form (3.10). We do not know how to apply this form and what is its generalization to the cases $n>2$.

\medskip
{\bf 4. Counterexample to $h^1(M)=h_1(M)$.}
\medskip
Let us consider the case $q=2$, $n=2$, $\ve=0$. We fix the following matrix $A=\left(\matrix \th & \th^6\\ \th^{-2}&0 \endmatrix \right)$.
The calculations below show that $h^1(M(A))=0$ (Theorem 4.1), $h^1(M(A^t))=1$ (Theorem 4.7). According to 1.10.3, this means that $h^1(M(A))=0$, $h_1(M(A))=1$ --- a counterexample to (0.3(a)).
\medskip
{\bf Theorem 4.1.} For the above $A$ we have $h^1(M(A))=0$.
\medskip
{\bf Proof.} For this $A$ the numbers $a_i$ of (3.9) are:

$$a_4=\th^{20}, \ \ \ a_3=\th^{16}, \ \ \ a_2=\th^{4}, \ \ \ a_1=\th^{3}, \ \ \ a_0=\th^{2};$$
$$\ord a_4=-20, \ \ \ \ord a_3=-16, \ \ \ \ord a_2=-4, \ \ \ \ord a_1=-3, \ \ \ \ord a_0=-2.$$
Later "equation (3.10)" will mean the equation (3.10) with these values of $a_i$ (and also values of $b_{ij}$ coming from $A$, see (4.4.1), (4.4.2) below). The Newton polygon for (3.10), $i=0$ has vertices $(1,-2)$; $(8,-16)$; $(16,-20)$. We denote elements of a $\n F_2$-basis of $S_0$ by $x_{j0}$, $j=1,...,4$, and solutions to (3.10) over them by $\{x_j\}=\sum_{i=0}^\infty x_{ji}T^i$. We have $\ord x_{j0}=2$, for $j=1,2,3$ and
$\ord x_{40}=\frac12$.

The equation (3.10) is not simple (see below), hence we need one more term for $x_{40}$.
\medskip
{\bf Lemma 4.2.}  $x_{40}=\th^{-\frac12}+\th^{-\frac{29}{16}}+\de_{41}$ where $\ord \de_{41}>\frac{29}{16}$ (the value of $\de_{41}$ depends on the value of $x_{40}$; all of them have ord $>\frac{29}{16}$).
\medskip
{\bf Proof.} First, we let
$$x_{40}=\th^{-\frac12}+\de_{40}\eqno{(4.2.1)}$$
where $\de_{40}$ is a new unknown. Substituting (4.2.1) to (3.10), $i=0$, we get
$$\sum_{j=0}^4 a_j\de_{40}^{2^j}+\sum_{j=0}^4 a_j(\th^{-\frac12})^{2^j}=0.\eqno{(4.2.2)}$$
We have $$\sum_{j=0}^4 a_j(\th^{-\frac12})^{2^j}=\th^{20}\th^{-8}+\th^{16}\th^{-4}+\th^{4}\th^{-2}+\th^{3}\th^{-1}+\th^{2}\th^{-\frac12}=\th^{\frac32}.$$
Hence, the Newton polygon of (4.2.2) has vertices $(0,-\frac32); \ (8,-16); \ (16, -20)$ and $\ord \de_{40}=\frac{29}{16}$ (8 values), $\frac12$ (8 values). There exist values of $\de_{40}$ such that $\th^{-\frac12}+\de_{40}=\sum_{j=1}^3 c_j x_{j0}$ where $c_j\in \n F_2$. For these $\de_{40}$ we have that their ord's are $\frac12$, exactly 8 values. Values of $\de_{40}$ such that $\ord \de_{40}=\frac{29}{16}$ give us 8 solutions to ((3.10), $i=0$) having ord $=\frac12$.

We can let
$$\de_{40}=\th^{-\frac{29}{16}}+\de_{41}.\eqno{(4.2.3)}$$
Substituting (4.2.3) to (4.2.2) we get
$$\sum_{j=0}^4 a_j\de_{41}^{2^j}+\sum_{j=0}^4 a_j(\th^{-\frac{29}{16}})^{2^j}+\th^{\frac32}=0.\eqno{(4.2.4)}$$
We have $$\sum_{j=0}^4 a_j(\th^{-\frac{29}{16}})^{2^j}+\th^{\frac32}= \th^{20}\th^{-29}+\th^{16}\th^{-\frac{29}{2}}+\th^{4}\th^{-\frac{29}{4}}+\th^{3}\th^{-\frac{29}{8}}+\th^{2}\th^{-\frac{29}{16}}+\th^{\frac32}= \th^{\frac3{16}}+ \delta $$ where $\ord \delta >0$.
Hence, the Newton polygon of (4.2.2) has vertices $(0,-\frac3{16})$; $ (8,-16)$; $ (16, -20)$ and $\ord \de_{41}=\frac{16-\frac3{16}}{8}$ (8 values), $\frac12$ (8 values). Case $\ord \de_{41}=\frac12$ gives us again values $\sum_{j=1}^3 c_j x_{j0}$ as earlier, case $\ord \de_{41}=\frac{16-\frac3{16}}{8}$ gives us the desired, because $\frac{16-\frac3{16}}{8}>\frac{29}{16}$. $\square$
\medskip
{\bf Remark 4.3.} Here we have a phenomenon observed by [C], see also [AB]: there are cases where the method of consecutive approximations does not give a solution to a polynomial. Namely, let us denote by ${\p}_{s}$ a subfield of $\p$ formed by series generated by rational powers of $\th^{-1}$, with coefficients in $\bar \n F_2$. More exactly, let $\al_1<\al_2<\al_3...$ be a sequence of rational numbers such that lim $\al_i=+\infty$, and $c_i\in \bar \n F_2$ coefficients. By definition, ${\p}_{s}$ is a proper subfield of $\p$ formed by all sums $\sum_{i=1}^\infty c_i\th^{-\al_i}$. A well-known example of $\g r\in \p-{\p}_{s}$ is a root to the Artin - Schreier polynomial $$x^2+x+\th^2=0\eqno{(4.3.0)}$$ (here $q=2$). Really, in a formal ring we have $$\g r=\th+\th^{\frac12}+\th^{\frac14}+\th^{\frac18}+...\eqno{(4.3.1)}$$ but this series $\not\in {\p}_{s}$. We let $\g r_i=\th+\th^{\frac12}+\th^{\frac14}+\th^{\frac18}+...+\th^{\frac1{2^n}}+\delta_{in}$, $i=1,2$ (there are two roots: (4.3.0) is separable). We have $\delta_{in}$ is a root to $$y^2+y+\th^{\frac1{2^n}}=0$$ and hence both $\delta_{1n}$, $\delta_{2n}$ have $\ord=-\frac1{2^{n+1}}$. This shows once again that the series (4.3.1) does not converge to $\g r$. It belongs to other field of power series, for example to $\bar \n F_p\langle \langle \th \rangle \rangle$ --- a field of power series with well-ordered support (see, for example, [P], 2.2.2).
\medskip
For $x_{40}$ we have exactly this phenomenon, but the first two terms of the equality $x_{40}=\th^{-\frac12}+\th^{-\frac{29}{16}}+\de_{41}$ are sufficient for our purpose.
\medskip
{\bf 4.4.} The same arguments as in Lemma 4.2 applied to $x_{j0}$, $j=1,2,3$, show that
$x_{j0}=c_j\th^{-2}+\de_j\in {\p}_{s}$ where $c_1,c_2,c_3$ form a basis of $\n F_8/\n F_2$ and $\ord \de_j>2$.
\medskip
Let us fix some $x_{40}$ and let us consider its minimal chain (see (2.5)). We have

$$b_{14}=\th^{16}+\th^{12}, \ \ b_{13}=\th^{12}, \ \ b_{12}=\th^{8}+\th^{2};\eqno{(4.4.1)}$$
$$\ \ord b_{14}=-16, \ \ \ord b_{13}=-12,\ \ \ord b_{12}=-8.$$

Hence, $$\sum_{k=2}^4b_{1k}x_{40}^{2^k}=(\th^{16}+\th^{12})(\th^{-\frac12}+\th^{-\frac{29}{16}}+\de_{41})^{16}+ \th^{12} (\th^{-\frac12}+\th^{-\frac{29}{16}}+\de_{41})^8 + $$ $$ +(\th^{8}+\th^{2})(\th^{-\frac12}+\th^{-\frac{29}{16}}+\de_{41})^4 = \th^6+ \delta$$ where $\ord \delta =-4$. This means that the Newton polygon of (3.10), $i=1$ has vertices  $(0,-6)$; $ (8,-16)$; $ (16, -20)$ and $\ord x_{41}=\frac54$ (the minimal value of $x_{41}$).

Finally, $$b_{24}=\th^8, \ \ord b_{24}=-8.\eqno{(4.4.2)}$$

Now we can use induction:
\medskip
{\bf Lemma 4.4a.} There exists a solution $\sum_{i=0}^\infty x_{4i}T^i$ to (3.10) having $\ord x_{4i}=2-\frac{3}{2^{i+1}}$.
\medskip
{\bf Proof.} We showed that this is true for $i=0,1$. Let us show that if this is true for $i=\al-2$ and $i=\al-1$ then this is true for $i=\al$. We have
$$\ord (b_{14}x_{4,\al-1}^{16})=16-\frac{48}{2^\al}; \ \ \ \ \ \ \ord (b_{13}x_{4,\al-1}^{8})=4-\frac{24}{2^\al};$$
$$\ord (b_{12}x_{4,\al-1}^{4})=-\frac{12}{2^\al};\ \ \ \ \ \ \ord (b_{24}x_{4,\al-2}^{16})=24-\frac{48}{2^{\al-1}}.$$
For $\al\ge2$ the minimal of these four numbers is the third one (really, for $\al=2$ this is true; for $\al\ge 3$ the third number is negative while the first, second and forth are positive), hence the Newton polygon of ((3.10), $i=\al$) has vertices  $(0,-\frac{12}{2^\al})$; $ (8,-16)$; $ (16, -20)$ and $\ord x_{4\al}=(16-\frac{12}{2^\al})/8=2-\frac{3}{2^{\al+1}}$. $\square$
\medskip
{\bf Lemma 4.5.} For $j=1,2,3$ there exist solutions $\{x_j\}$ over $x_{j0}$ having $\forall \ i$ \ $\ord x_{ji}=2$.
\medskip
{\bf Proof.} For $j=1,2,3$ we have $\ord x_{j0}=2$, hence for these $j$
$$\ord (b_{14}x_{j0}^{16})=16;\ \ \ \ord (b_{13}x_{j0}^{8})=4;\ \ \ \ord (b_{12}x_{j0}^{4})=0$$
and the Newton polygon of (3.10), $i=1$, $j=1,2,3$ has vertices  $(0,0)$; $ (8,-16)$; $ (16, -20)$ and $\ord x_{j1}$ can be chosen 2.

The same situation holds for $i\ge2$. We have
$$\ord (b_{24}x_{j0}^{16})=24,$$
hence by induction we get that there are solutions $\sum_{i=0}^\infty x_{ji}T^i$ to (3.10), $j=1,2,3$, having $\ord x_{ji}=2$. $\square$
\medskip
So, we got 4 basis solutions $\{x_j\}=\sum_{i=0}^\infty x_{ji}T^i$, $j=1,...,4$. Any solution to (3.10) is $\sum_{j=1}^4 C_j \{x_j\}$ where $C_j\in \n F_2[[T]]$ (Proposition 2.3).
\medskip
{\bf Lemma 4.6.} The set $\sum_{j=1}^4 C_j \{x_j\}$ does not contain small solutions (here and below --- except the zero solution).
\medskip
{\bf Proof.} Let us assume that $\exists \ C_1,...,C_4$ such that $\sum_{j=1}^4 C_j \{x_j\}$ is a small solution. We consider $S_{123}:=\sum_{j=1}^3 C_j \{x_j\}$, we denote $S_{123}=\sum_{i=0}^\infty \bar x_{1,2,3;i}T^i$. We have: $\ord \bar x_{1,2,3;i}\ge2$, because $\forall \ i$ elements $\bar x_{1,2,3;i}$ are linear combinations of $x_{jk}$ for $j=1,2,3$, $k\le i$ with coefficients in $\n F_2$.

Further, we denote $S_{4}:=C_4 \{x_4\}=\sum_{i=0}^\infty \bar x_{4i}T^i$. Lemma 4.4a shows that $\forall \ i$ $\ord x_{4i}$ are different and $\frac12\le \ord x_{4i}<2$, hence $\forall \ i$ we have $\frac12\le \ord \bar x_{4i}<2$. This means that $\sum_{j=1}^4 C_j \{x_j\}=S_{123}+S_4$ cannot be a small solution. $\square$
\medskip
This gives us a proof of Theorem 4.1. Really, if $H^1(M(A))\ne0$ then equation (3.10) has a small solution. $\square$
\medskip
{\bf Remark.} It is easy to show that if $\bar x_{1,2,3;0} \ne0$ then $\forall \ i>0$ we have $\ord \bar x_{1,2,3;i}=2$. Really, let us assume that $\exists \ i $ such that $\ord \bar x_{1,2,3;i}>2$. We choose minimal such $i$, and denote it by $i_0$. The condition $\bar x_{1,2,3;0} \ne0$ implies $i_0\ge1$. The calculation of Lemma 4.5 shows that --- because for $i=0,...,i_0-1$ we have $\ord \bar x_{1,2,3;i}=2$, we have either $\ord \bar x_{1,2,3;i_0}=2$ or $\ord \bar x_{1,2,3;i_0}=\frac12$. The condition $\ord \bar x_{1,2,3;i_0}=\frac12$ contradicts to $\ord \bar x_{1,2,3;i}\ge2$, and the condition $\ord \bar x_{1,2,3;i_0}=2$ contradicts to the choice of $i_0$.
\medskip
{\bf Theorem 4.7.} For the $A$ above we have $h^1(M(A^t))=1$.
\medskip
{\bf Proof.} We have $A^t=\left(\matrix \th & \th^{-2}\\ \th^{6}&0 \endmatrix \right)$. The numbers $a_i$ of (3.10) for $A^t$ are:

$$a_4=\th^{-12}, \ \ \ a_3=\th^{-16}, \ \ \ a_2=\th^{-20}, \ \ \ a_1=\th^{-5}, \ \ \ a_0=\th^{-6};$$
$$\ord a_4=12, \ \ \ \ord a_3=16, \ \ \ \ord a_2=20, \ \ \ \ord a_1=5, \ \ \ \ord a_0=6.$$
Now "equation (3.10)" will mean the equation (3.10) with these values of $a_i$ and the below values of $b_{ij}$.
The Newton polygon for (3.10), $i=0$ has vertices (1,6); (2,5); (16,12), hence $\ord x_{10}=1$, and for $j=2,3,4$
we have $\ord x_{j0}=-\frac12$. We have:

$$b_{14}=\th^{-16}+\th^{-20}, \ \ b_{13}=\th^{-20}, \ \ b_{12}=\th^{-6}+\th^{-24},\ \ b_{24}=\th^{-24};$$
$$\ \ord b_{14}=16, \ \ \ord b_{13}=20,\ \ \ord b_{12}=6,\ \ \ord b_{24}=24,$$
hence the Newton polygon for (3.10), $i=1$, $x_0=x_{10}$ has vertices (0,10); (1,6); (2,5); (16,12), and we can choose $x_{11}$ having $\ord=4$.
\medskip
{\bf Lemma 4.8.} For $\{x_1\}$ we have: $\ord x_{1n}=4^n$.
\medskip
{\bf Proof.} Induction. Let the lemma hold for $n = i$ and $n = i + 1$. Then it holds for $n = i + 2$. Really,
$$\ord b_{14}x_{1,i+1}^{16}=16+16\cdot4^{i+1}; \ \ \ \ \ \ord b_{13}x_{1,i+1}^{8}=20+8\cdot4^{i+1};$$
$$\ord b_{12}x_{1,i+1}^{4}=6+4\cdot4^{i+1}; \ \ \ \ \ \ord b_{24}x_{1,i}^{16}=24+16\cdot4^{i}.$$

The minimal of these four numbers is the third one, it is $6+4^{i+2}$, this is the $y$-coordinate of the vertex of the Newton polygon having $x=0$. Hence, $\ord x_{i+2}=6+4^{i+2}-6=4^{i+2}$. $\square$
\medskip
{\bf Lemma 4.9.} Let $\{x\} = \sum_{i=0}^\infty x_{i}T^i$ be any solution to (3.10) over $x_{0}$ such that $\ord x_0=-\frac12$. Then $\forall \ i>0$ we have: $\ord x_{i}=-\frac12$.
\medskip
{\bf Proof.} Immediate, by induction. In notations of Lemma 4.8, we have
$$\ord b_{14}x_{i+1}^{16}=8; \ \ \ \ \ \ord b_{13}x_{i+1}^{8}=16;$$
$$\ord b_{12}x_{i+1}^{4}=4; \ \ \ \ \ \ord b_{24}x_{i}^{16}=16,$$
hence the Newton polygon of (3.10) for any $i$ is one segment $(0,4) - (16,12)$, hence the lemma. $\square$
\medskip
{\bf Lemma 4.10.} Any small solution to (3.10) belongs to $\n F_2[T]\{x_1\}$.
\medskip
{\bf Proof.} Let $y=(y_0,y_1,\dots)$ be a small solution. Lemma 4.9 implies that $y_0\in \n F_2 \ x_{10}$, i.e. $\exists \ k_0\in \n F_2$ such that $y_0=k_0 x_{10}$. Let us consider $y-k_0\{x_1\}$. It is a small solution, its first term is 0, and hence we can divide it by $T$. Now we continue the process: there exists $k_1\in \n F_2$ such that the first term of $(y-k_0\{x_1\})/T$ is $k_1 x_{10}$. We consider $((y-k_0\{x_1\})/T-k_1\{x_1\})/T$ etc. As a result, we get that $\exists K:=\sum_{n=0}^\infty  k_nT^n\in \n F_2[[T]]$ such that $y= \{x_1\}K$. It is easy to see that $K\in \n F_2[T]$. Really, let $k_n\ne0$. In this case the ord of the $n$-th term of $\{x_1\}K$ is 1 ( = $\ord x_{10}$), because $\forall \ i>0$ we have $\ord x_{10}>\ord x_{1i}$. The condition that $y$ is small implies that there exists only finitely many $k_n$ such that $k_n\ne0$. $\square$
\medskip
{\bf Remark.} There is another proof of this lemma. In the notations of Lemma 4.6, let $\sum_{j=1}^4 C_j \{x_j\}$ be a small solution. We denote $C_j=\sum_{i=0}^\infty c_{ji}T^i$, where $c_{ji}\in \n F_2$. Like in Lemma 4.6, we denote $S_{234}:=\sum_{j=2}^4 C_j \{x_j\}=\sum_{i=0}^\infty \bar x_{2,3,4;i}T^i$ and $S_{1}:= C_1 \{x_1\}=\sum_{i=0}^\infty \bar x_{1i}T^i$. If $S_{234}\ne0$ then there exists the minimal $i_0$ such that $\bar x_{2,3,4;i_0}\ne0$. This implies $\forall \ j=2,3,4, \forall \ k<i_0$ $c_{jk}=0$ and $\exists \ j=2,3,4$ such that $c_{ji_0}\ne0$. This means that $\ord \bar x_{2,3,4;i_0}=-\frac12$, and hence, according Lemma 4.9, $\forall \ i\ge i_0$ we have $\ord \bar x_{2,3,4;i}=-\frac12$. Like in the proof of Lemma 4.6, we get that $\forall \ i$ $\ord \bar x_{1i}\ge1$, hence the sum $S_{1}+S_{234}$ cannot be a small solution. The only exception is $S_{234}=0$. If $C_1\not\in \n F_2[T]$ then there exists infinitely many $i$ such that $\ord \bar x_{1i}=1$ (again because all $\ord x_{1i}$ are different), hence for a small solution $C_1$ must belong to $\n F_2[T]$.
\medskip
{\bf 4.11. End of the proof.} Lemma 4.10 implies that $h^1(M(A^t))=1$. Really, (3.10) for the present case means that the set of small $y_{21}$ of Section 3 has dimension 1. Further, (3.4) shows that if $y_{21}$ is small then $y_{22}$ is also small, and (3.0) shows that $y_1=(y_{11}, y_{12})$ are also small. $\square$
\medskip
%\newpage
{\bf 5. A question of D. Goss.}
\medskip
Prof. David Goss wrote ([G1])
\medskip
{\it One last question:  Let $\phi$ be a t-motive where $\phi_T$ has infinitesimal part $T+N$ where $N$ is unipotent. Define $\hat \phi$  to be generated by $\hat \phi_T$ where
$\hat \phi_T$ has exactly the same coefficients as $\phi_T$ but where $N$ is now set to $0$.
What is the relationship between these two objects? if one is uniformizable what about the other? And then what would be the relationship between the lattices? etc.
\medskip
David}
\medskip
We give an example that two t-motives mentioned above can be of different uniformizability type. We use notations of (1.10.1), $n=2$, $N$ of Section 3, i.e. $N=\ve \left(\matrix 0&1\\ 0&0 \endmatrix \right)$ where $\ve=0$ or 1. We consider the case $q>2$, $A:=\left(\matrix 0&a_{12}\\ a_{21}&0 \endmatrix \right)$ where $a_{21}$ is any number having ord $a_{21}=-\frac{q^2}{q-1}$, and $a_{12}$ satisfies $$\frac{\th^q}{a_{21}}+\frac{\th^{q^2}}{a_{21}^{q^2}}-a_{12}^q=0$$ (this expression is $a_2$ of (3.9)).
\medskip
{\bf Theorem 5.1.} For $\ve=0$ the t-motive $M(A,N)$ is non-uniformizable, while for $\ve=1$ it is uniformizable.
\medskip
{\bf Proof.} For $\ve=0$, resp. 1 we denote the corresponding $M(A,N)$ by $M_0$, resp. $M_1$. The coefficients (3.9) become
$$a_4=\frac{\th^{q^3+q^2}}{a_{21}^{q^2}};\ \ \ a_3=-\ve; \ \ \ a_2=a_1=0;\ \ \ a_0=\frac1{a_{21}};$$
$$\ord a_4=\frac{q^2}{q-1};\ \ \  \ord a_3=0\ (\ve=1);\ \ord a_3=+\infty\ (\ve=0);\ \ \  \ord a_0=\frac{q^2}{q-1};$$
$$b_{14}=-\frac{\th^{q^3}+\th^{q^2}}{a_{21}^{q^2}}; \ \ \ b_{13}=0; \ \ \  b_{12}=\frac{1}{a_{21}}+\frac{1}{a_{21}^{q^2}}; \ \ \ b_{24}=\frac{1}{a_{21}^{q^2}};$$
$$\ord b_{14}=\frac{q^2}{q-1}+q^2; \ \ \ \ord b_{12}=\frac{q^2}{q-1}; \ \ \ \ord b_{24}=\frac{q^4}{q-1}.$$

Let us show $M_0$ is not uniformizable. The Newton polygon of ((3.10), $i=0$) is a segment whose ends have coordinates
\medskip
$(1, \frac{q^2}{q-1})$, $(q^4, \frac{q^2}{q-1})$,
\medskip
hence all $x_0\ne0$ have $\ord=0$. We get by induction by $i$ that for any $i_0>0$ the Newton polygon for the equation ((3.10), $i=i_0$) is a segment whose ends have coordinates
\medskip
$(0, \frac{q^2}{q-1})$, $(q^4, \frac{q^2}{q-1})$,
\medskip
and hence (if $x_0\ne0$) for all $i $ we have ord $x_i=0$. This means that $H^1(M_0)=0$.
\medskip
Now let us show $M_1$ is uniformizable. The vertices of the Newton polygon for the equation ((3.10), $i=0$) have coordinates
\medskip
$(1, \frac{q^2}{q-1})$, $(q^3,0)$, $(q^4, \frac{q^2}{q-1})$,
\medskip
and hence the set of solutions to ((3.10), $i=0$) is a $\n F_q$-vector space of dimension 4 having a basis $x_{10},...,x_{40}$ such $\ord x_{j0}=\frac{q^2}{(q-1)(q^3-1)}$ for $j\le3$, $\ord x_{40}=-\frac{1}{(q-1)(q^2-q)}$.
\medskip
Let us consider the equation ((3.10), $i=1$) for $x_0=x_{40}$ (the "worst" case). We have:
\medskip
$\ord b_{14} x_{40}^{q^4}=\frac{q^4-2q^3}{(q-1)^2}; \ \ \ \ord b_{12} x_{40}^{q^2}=\frac{q^3-q^2-q}{(q-1)^2}$.
\medskip
Hence, for $q>2$ the vertices of the Newton polygon for the equation ((3.10), $i=1$), $x_0=x_{40}$ have coordinates
\medskip
$(0, \frac{q^3-q^2-q}{(q-1)^2})$, $(q^3,0)$, $(q^4, \frac{q^2}{q-1})$,
\medskip
and hence ((3.10), $i=1$) has a solution $x_{41}$ having $\ord=\frac{q^2-q-1}{(q-1)^2q^2}$.

Hence, we have
\medskip
$\ord b_{14} x_{41}^{q^4}=\frac{2q^4-2q^3-q^2}{(q-1)^2}; \ \ \ \ord b_{12}x_{41}^{q^2}=\frac{q^3-q-1}{(q-1)^2};\ \ \ \ord b_{24}x_{40}^{q^4}=\frac{q^5-q^4-q^3}{(q-1)^2}.$
\medskip
For $q>2$ the minimal of these 3 numbers is $\frac{q^3-q-1}{(q-1)^2}$, hence the vertices of the Newton polygon for ((3.10), $i=2$) have coordinates
\medskip
$(0,\frac{q^3-q-1}{(q-1)^2})$, $(1,\frac{q^2}{q-1})$, $(q^3,0)$, $(q^4,\frac{q^2}{q-1})$,
\medskip
and hence ((3.10), $i=2$) has a solution $x_{42}$ having $\ord=\frac{q^2-q-1}{(q-1)^2}$. This is sufficient to use induction. Namely, let us prove that
\medskip
(*) $\forall n\ge1$ we have $\ord x_{4n}\ge \frac12 q^{2(n-2)}$.
\medskip
(*) holds for $n=1,2$, hence to prove (*) for all $n$ it is sufficient to prove:
\medskip
{\bf Lemma 5.2.} Let (*) hold for $n=i$ and $n=i+1$. Then it holds for $n=i+2$.
\medskip
{\bf Proof.} It is straightforward. We have
\medskip
$\ord b_{14}x_{4,i+1}^{q^4}\ge \frac{q^3}{q-1}\cdot\frac12\cdot q^{2(i-1)}\cdot q^4> \frac12 q^{2i+4}$;
\medskip
$\ord b_{12}x_{4,i+1}^{q^2}\ge \frac{q^2}{q-1}\cdot\frac12\cdot q^{2(i-1)}\cdot q^2> \frac12 q^{2i+1}$;
\medskip
$\ord b_{24}x_{4i}^{q^4}\ge \frac{q^4}{q-1}\cdot\frac12\cdot q^{2(i-2)}\cdot q^4> \frac12 q^{2i+3}$.
\medskip
We have $\frac12 q^{2i+1}-\frac{q^2}{q-1}>\frac12 q^{2i}$, hence the lemma. $\square$
\medskip
Calculations for minimal chains of $x_{j0}$, $j\le3$, are similar: the $y$-coordinate of the leftmost vertex of the Newton polygon for (3.10), any $i$ is greater than the same coordinate for the above minimal chain of $x_{40}$. Alternatively, we can prove that the above minimal chain of $x_{40}$ is simple, and to apply Proposition 2.7. The details are left to the reader.
\medskip
Finally, the same arguments as in (4.11) show that if we have 4 linearly independent small solutions $y_{21}$ to (3.10), then they give 4 linearly independent small solutions to (3.0a). This means that $M_1$ is uniformizable. $\square$
\medskip
{\bf 6. Appendix.}
\medskip
{\bf 6.1.} Here we give an analog of calculations of Section 3 for a non-pure t-motive of dimension 2 and rank 5. In terminology of [GL07], it is a standard-1 t-motive from ([GL07], 11.1) having $\la_1=3, \ \la_2=2$. (1.3) for it is the following:
$$T\left(\matrix e_1\\ e_2 \endmatrix \right)=\th \left(\matrix e_1\\ e_2 \endmatrix \right) + \left(\matrix a_{11} & a_{12} \\ a_{21} & a_{22} \endmatrix \right)\tau \left(\matrix e_1\\ e_2 \endmatrix \right) + \left(\matrix b_{1} & 0 \\ b_{2} & 1 \endmatrix \right)\tau^2\left(\matrix e_1\\ e_2 \endmatrix \right) + \left(\matrix 1 & 0 \\ 0& 0 \endmatrix \right) \tau^3\left(\matrix e_1\\ e_2 \endmatrix \right).$$ Its $f_*$-basis can be chosen as

$\left(\matrix e_1 \\ e_{2} \\ \tau e_1 \\ \tau e_2 \\ \tau^2e_1 \endmatrix \right)$. The matrix $Q$ in this basis is $\left(\matrix \\ 0_{32} & I_3 \\ \\ {\matrix 0 & T-\th \\ T-\th & 0\endmatrix} & {\matrix -a_{21}& -a_{22}& -b_{2}\\ -a_{11}& -a_{12}& -b_{1} \endmatrix} \endmatrix \right)$ where $0_{32}$, $ I_3$ are respectively $3\times2$, $3\times3$-blocks. We have (here $Y=(y_1,\dots,y_5)$ ):
$$y_5^{(1)}=-\frac1{a_{12}}y_4 -\frac{a_{22}}{a_{12}}y_4^{(1)} +\frac{T-\th^q}{a_{12}}y_4^{(2)};\eqno{(6.1.1)}$$
$$y_5+b_1y_5^{(1)}+a_{11}^qy_5^{(2)}-(T-\th^{q^2})y_5^{(3)}+b_2y_4^{(1)}+a_{21}^qy_4^{(2)}=0.\eqno{(6.2.1)}$$
Further, we transform (here $(6.\al.*)$ come from $(6.\al.1)$, $\al=1,\ 2$):
$$y_5^{(1)}+b_1^qy_5^{(2)}+a_{11}^{q^2}y_5^{(3)}-(T-\th^{q^3})y_5^{(4)}+b_2^qy_4^{(2)}+a_{21}^{q^2}y_4^{(3)}=0;\eqno{(6.2.2)}$$
$$y_5^{(2)}=-\frac1{a_{12}^{q}}y_4^{(1)} -\frac{a_{22}^{q}}{a_{12}^{q}}y_4^{(2)} +\frac{T-\th^{q^2}}{a_{12}^{q}}y_4^{(3)};\eqno{(6.1.2)}$$
$$y_5^{(3)}=-\frac1{a_{12}^{q^2}}y_4^{(2)} -\frac{a_{22}^{q^2}}{a_{12}^{q^2}}y_4^{(3)} +\frac{T-\th^{q^3}}{a_{12}^{q^2}}y_4^{(4)};\eqno{(6.1.3)}$$
$$y_5^{(4)}=-\frac1{a_{12}^{q^3}}y_4^{(3)} -\frac{a_{22}^{q^3}}{a_{12}^{q^3}}y_4^{(4)} +\frac{T-\th^{q^4}}{a_{12}^{q^3}}y_4^{(5)}.\eqno{(6.1.4)}$$
Now we substitute (6.1.1) -- (6.1.4) to (6.2.2) in order to eliminate $y_5$:
$$ -\frac1{a_{12}}y_4 -\frac{a_{22}}{a_{12}}y_4^{(1)} +\frac{T-\th^q}{a_{12}}y_4^{(2)} +b_1^q  (-\frac1{a_{12}^{q}}y_4^{(1)} -\frac{a_{22}^{q}}{a_{12}^{q}}y_4^{(2)} +\frac{T-\th^{q^2}}{a_{12}^{q}}y_4^{(3)}) $$ $$+a_{11}^{q^2} (-\frac1{a_{12}^{q^2}}y_4^{(2)} -\frac{a_{22}^{q^2}}{a_{12}^{q^2}}y_4^{(3)} +\frac{T-\th^{q^3}}{a_{12}^{q^2}}y_4^{(4)}) -(T-\th^{q^3}) (-\frac1{a_{12}^{q^3}}y_4^{(3)} -\frac{a_{22}^{q^3}}{a_{12}^{q^3}}y_4^{(4)} +\frac{T-\th^{q^4}}{a_{12}^{q^3}}y_4^{(5)}) $$ $$ +b_2^qy_4^{(2)}+a_{21}^{q^2}y_4^{(3)}=0.\eqno{(6.2.3)}$$
We get an equation of type (2.2.2) having $r=5$, $n=2$, that supports Conjecture 2.3.6.
\medskip
{\bf Remark 6.3.} For t-motives with higher $r$, $n$ the similar calculation gives us systems of type

$$\sum_{\ga=0}^{r_{11}} a_{1\ga} x^{(\ga)} + \sum_{\ga=0}^{r_{12}} b_{1\ga} y^{(\ga)}=0;\eqno{(6.3.1)}$$
$$\sum_{\ga=0}^{r_{21}} a_{2\ga} x^{(\ga)} + \sum_{\ga=0}^{r_{22}} b_{2\ga} y^{(\ga)}=0.\eqno{(6.3.2)}$$ where $a_{ij}$, $b_{ij}\in \p(T)$, $x, \ y\in \p[[T]]$ are unknowns. For elimination of an unknown from this system we can use the theory of the $p$-resultant, see, for example, [G], Section 1.5. We expect that for any explicitly given t-motive (for example for standard-2 t-motives of [GL07], 11.2) the similar calculations will give us a proof of Conjecture 2.13 for them.
\medskip
{\bf Example 6.4.} Let us consider $M$ from (1.10.1) for $n=3$, $N=0$. Notations and calculations are similar to the ones of Section 3. We write $Y=(y_1, y_2)$ as a block matrix where $y_1=(y_{11}, y_{12}, y_{13})$, $y_2=( y_{21}, y_{22}, y_{23})$. Analog of (3.1) is the same, (3.2) and (3.3) become (we give here a sketch of calculations):
$$y_{21}=y_{21}^{(2)}(T-\th^q)-y_{21}^{(1)}a_{11}-y_{22}^{(1)}a_{21}-y_{23}^{(1)}a_{31};\eqno{(6.4.1)}$$
$$y_{22}=y_{22}^{(2)}(T-\th^q)-y_{21}^{(1)}a_{12}-y_{22}^{(1)}a_{22}-y_{23}^{(1)}a_{32};\eqno{(6.4.2)}$$
$$y_{23}=y_{23}^{(2)}(T-\th^q)-y_{21}^{(1)}a_{13}-y_{22}^{(1)}a_{23}-y_{23}^{(1)}a_{33}.\eqno{(6.4.3)}$$
Analogs of (3.4) - (3.6) become
$$y_{23}^{(1)}=-\frac{1}{a_{31}}y_{21}+\frac{T-\th^q}{a_{31}}y_{21}^{(2)}-\frac{a_{11}}{a_{31}}y_{21}^{(1)} -\frac{a_{21}}{a_{31}}y_{22}^{(1)};\eqno{(6.4.4)}$$
$$y_{23}^{(2)}=-\frac{1}{a_{31}^q}y_{21}^{(1)}+\frac{T-\th^{q^2}}{a_{31}^q}y_{21}^{(3)}-\frac{a_{11}^q}{a_{31}^q}y_{21}^{(2)} -\frac{a_{21}^q}{a_{31}^q}y_{22}^{(2)};\eqno{(6.4.5)}$$
$$y_{23}^{(3)}=-\frac{1}{a_{31}^{q^2}}y_{21}^{(2)}+\frac{T-\th^{q^3}}{a_{31}^{q^2}}y_{21}^{(4)}-\frac{a_{11}^{q^2}}{a_{31}^{q^2}}y_{21}^{(3)} -\frac{a_{21}^{q^2}}{a_{31}^{q^2}}y_{22}^{(3)}.\eqno{(6.4.6)}$$
Applying $\tau$ to (6.4.3) we get an analog of (3.7):
$$y_{23}^{(1)}=y_{23}^{(3)}(T-\th^{q^2})-y_{21}^{(2)}a_{13}^q-y_{22}^{(2)}a_{23}^q-y_{23}^{(2)}a_{33}^q.\eqno{(6.4.7)}$$

Substituting (6.4.4) - (6.4.6) to (6.4.2) and (6.4.7), we eliminate $y_{23}$ from (6.4.1) - (6.4.3). We get a system of type (6.3.1) - (6.3.2) with unknowns $y_{21}$, $y_{22}$, having $r_{11}, r_{12}, r_{21}, r_{22}$ respectively 2, 2, 4, 3 (and $b_{20}=0$). We can eliminate one of these two unknowns either using the theory of the $p$-resultant, or we can enlarge this system (6.3.1) - (6.3.2) applying $\tau$, $\tau^2$, $\tau^3$ to (6.3.1) and $\tau$ to (6.3.2). We shall get a matrix equation $$\Cal A\left(\matrix y_{21}\\ \\{y_{21}}^{(1)}\\ \dots\\ {y_{21}}^{(5)} \endmatrix \right) +\Cal  B\left(\matrix y_{22}\\ \\{y_{22}}^{(1)}\\ \dots\\ {y_{22}}^{(5)} \endmatrix \right)=0\eqno{(6.4.8)}$$ where $\Cal A$, $\Cal B\in M_{6\times 6}(\p[T])$. We have $\Cal A\in GL_6(\p(T))$; multiplying (6.4.8) from the left by the first line of $\Cal A^{-1}$ we get an expression of $y_{21}$ as a linear combination of $y_{22}, y_{22}^{(1)}, \dots, y_{22}^{(5)}$. Substituting this expression to (6.3.1) we get an equation
$$\sum_{\ga=0}^{7} c_{\ga} y_{21}^{(\ga)}=0\eqno{(6.4.9)}$$ for the unknown $y_{21}$, where $c_\ga\in \p[T]$.
\medskip
Why the maximal value of $\ga$ is 7 and not 6, and finding of the value of $n$ (see Conjecture 2.3.6) is a subject of further calculations.
\medskip
{\bf 6.5. General case.} Here we consider constructions of (1.4.3.1), (1.4.3.2), (1.5.1) in some general setting. Let $R=R_1\subset S=R_2\subset \g S=R_3$ be commutative rings, and $\si$ concordant automorphisms of $R_i$.
\medskip
{\bf 6.5.1.} In (1.4.3.1), (1.4.3.2), (1.5.1) we have $R=\p[T], \ S=\p\{T\}, \ \g S=\p[[T]]$, and for $z\in \p[[T]]$ we have $\si(z)=z^{(1)}$.
\medskip
We consider non-commutative rings $R_i\{\tau\}$ having the commutation rule $\tau r=\si(r)\cdot\tau$. Let $M$ be a $R$-module with $\tau$-action such that $\tau$ is an automorphism of $M$ and such that $$\tau(rm)=\si(r)\cdot\tau(m)\eqno{(6.5.2)}$$ holds. This condition implies that we can consider $M$ as a $R\{\tau\}$-module. We shall consider only these modules.
\medskip
Now we denote $M^S:=M\otimes_RS$, $M_S:=\Hom_R(M,S)$. They have a natural structure of $S$-modules. $\tau$ acts on them by the standard formulas $$\tau(m\otimes s)=\tau(m)\otimes \si(s), \ \ \ (\tau(\vf))(m)=\si(\vf(\tau^{-1}(m)))$$ It is immediately checked that condition (6.5.2) holds for them (where $r\in S$, $m\in M^S$ or $m\in M_S$), hence $M^S$, $M_S$ are $S\{\tau\}$-modules.
\medskip
For a $\Cal R\{\tau\}$-module $\Cal M$ we have $\Cal M^\tau=\Hom_{\Cal R\{\tau\}}(\Cal R,\Cal M)$. Particularly (we let $\Cal R=S$, $\Cal M=M_S$), $${M_S}^\tau=\Hom_{S\{\tau\}}(S,\Hom_R(M,S))=\Hom_{R\{\tau\}}(M,S)$$
\medskip
{\bf Lemma 6.5.3}. $\forall \ i$ we have $\Ext_{S\{\tau\}}^i(S,\Hom_R(M,S))=\Ext_{R\{\tau\}}^i(M,S)$.
\medskip
{\bf Proof} is an exercise for the reader. $\square$
\medskip
Analogously, we can consider groups $\Ext_{S\{\tau\}}^i(S,M^S)$. Are they interesting objects for our $R$, $S$, $M$?
\medskip
Further, for a $\Cal R\{\tau\}$-module $\Cal M$ we have $$\Cal M_\tau:=\Cal R\underset{\Cal R\{\tau\}}\to{\otimes}\Cal M=\Cal M/<\tau(m)-m>$$ where $<\tau(m)-m>$ is a submodule of $\Cal M$ generated by all these elements ($m\in \Cal M$).
\medskip
For $M=\Cal M$ uniformizable we have: both ${M^S}_\tau$, ${M_S}_\tau$ are 0. Maybe for other $M$ they (or the corresponding Ext's) are non-zero?
\medskip

{\bf 6.6. Further research.} We want to answer Question 0.3c. Explicitly, what is the set of quadruples $(c_1,\dots,c_4)$ such that there exists a (pure, having $N=0$) t-motive $M$ such that $c_1=r(M)$, $c_2=h^1(M)$, $c_3=h_1(M)$ and $c_4$ is the rank of the pairing $\pi$? Trivial restrictions are: for uniformizable $M$ all $c_i$ must be equal, for non-uniformizable $M$ all $c_i$ must satisfy $c_1>c_2$, $c_1>c_3$, $c_2\ge c_4\ge0$, $c_3\ge c_4$, and $$c_4\ge c_2+c_3-c_1\eqno{(6.6.1)}$$ (see Section (6.7) below for the justification of (6.6.1)). Maybe these are the only restrictions on $c_*$? (Exception is the trivial case $c_1\le3$: really, all pure, having $N=0$ t-motives of rank $\le3$ are either Drinfeld modules or their duals, hence they are all uniformizable).
\medskip
For the case $r=c_1=4$ we can try to find all quadruples satisfying $4> c_2,c_3\ge c_4\ge0$ and (6.6.1) by computer search among t-motives of the form (1.10.1), $N=0$, $n=2$. For higher $r$ we can either consider $M$ from (6.1) (consideration of t-motives of the form (1.10.1), $N=0$, $n=3$ can be too complicated, see 6.4), or to consider direct sums of the above $M$ and some Drinfeld modules. Since all invariants $r(M)$, $h^1(M)$, $h_1(M)$ and the rank of the pairing $\pi$ are additive with respect to the direct sum, most likely it will be enough to find finitely many small quadruples in order to get all possible quadruples as their sums.
\medskip
{\bf 6.6.2.} Also, we can calculate dim $\Ext^1_{\p[T,\tau]}(M,Z_1)$ (notations of [G], (5.9.22); see [G], Remark 5.9.26 for the meaning of this space). If $M$ is uniformizable then this Ext is 0. Also, there are some other Ext modules (see 6.5) related to $M$, as well as Tor groups, for example $\Tor_1^{\p[T,\tau]}(M,Z_1)$. Are there relations between $h^1(M)$, $h_1(M)$ and the dimensions of these Ext, Tor?
\medskip
There is also a technical problem to prove or disprove Conjecture 2.13. As a first step, we should understand why the equation (6.4.9) has degree 7 and not 6 as it should be. What is the dimension of the set of solutions to (6.3.1) -- (6.3.2)? What is a meaning (and a generalization) of equation (3.11)?
\medskip
Further, we give in Section 5 an example of $M_0$, $M_1$ given by the formula (1.10.1) such that the matrices $A$ for $M_0$, $M_1$ coincide, $N$ for $M_0$ is 0 while $N$ for $M_1$ is not 0, and such that $M_0$ is non-uniformizable, while $M_1$ is uniformizable. Are there examples $M$ such that the situation is inverse, i.e. $M_0$ is uniformizable and $M_1$ is non-uniformizable?
\medskip
Finally, we can try to find $h^1$, $h_1$ for all $M$ described by (1.10.1), first for $n=2$, and to describe explicitly the set of all $2\times2$-matrices $A$ such that $M(A)$ from (1.10.1) is uniformizable. For example, what is the minimal value of $u$ having the property:
\medskip
If ord of all entries of $A$ is $>u$ (or $\ge u$) then $M(A)$ is uniformizable?
\medskip
[GL17], Proposition 2 shows that $u\ge\frac{q}{q^2-1}$; most likely this bound can be improved. Having a list of all such $A$ we can try to find explicitly Siegel matrices of the lattices of $M(A)$ and to check whether all lattices of dimension 4 in $\p^2$ can be obtained by this manner or not. This is a next step to a problem whether the lattice map of pure uniformizable t-motives is surjective (or near-surjective), or not (see [GL17], Introduction for a discussion of the isomorphism problem).
\medskip
{\bf 6.7. Justification of (6.6.1).} Let us consider equations (1.6.1), (1.6.2) for the case $X\in M_{r\times1}(\p[[T]])$, $Y\in M_{1\times r}(\p[[T]])$, we denote them by $(1.6.1_\infty)$, $(1.6.2_\infty)$ respectively. They coincide with (1.7.1) for $V=M[[T]]$ and the dual module respectively. The set of their solutions is a $\n F_q[[T]]$-module. (1.7.3), (1.7.4) define a $\n F_q[[T]]$-pairing between them; we denote it by $\pi_\infty$.
\medskip
{\bf Proposition 6.7.1.} The set of solutions to $(1.6.1_\infty)$, $(1.6.2_\infty)$ is a free $\n F_q[[T]]$-module of dimension $r$. The pairing $\pi_\infty$ is perfect over $\n F_q[[T]]$.
\medskip
{\bf Proof.} Let us consider equations $(1.6.1_\infty)$, $(1.6.2_\infty)$ in the matrix form, i.e. let $\Cal Y, \ \Cal X\in M_{r\times r}(\p[[T]])$ be unknowns and
$$\Cal Y^{(1)}Q=\Cal Y, \eqno{(1.6.1m)}$$ $$Q\Cal X=\Cal X^{(1)} {(1.6.2m)}$$ equations (here and below $m$ means a matrix form). Practically, $\Cal Y$ is $\g Y$ from 1.7.5.1; we use another notation because we assume in 1.7.5.1 that $h^1=r$; here we consider a general case.
\medskip
Any line of $\Cal Y$ (resp. column of $\Cal X$) satisfying $(1.6.1m), \ (1.6.2m)$ is a solution to $(1.6.1_\infty)$, $(1.6.2_\infty)$. We denote $Q=Q_0+Q_1T+...+Q_\vk T^\vk$,
\medskip
$X=X_{0}+X_{1}T+X_{2}T^2+...$, \ $Y=Y_{0}+Y_{1}T+Y_{2}T^2+...$,
\medskip
$\Cal X=\Cal X_{0}+\Cal X_{1}T+\Cal X_{2}T^2+...$, \ $\Cal Y=\Cal Y_{0}+\Cal Y_{1}T+\Cal Y_{2}T^2+...$ where entries of $Q_i$, $X_{i}$, $Y_{i}$, $\Cal X_{i}$, $\Cal Y_{i}$ belong to $\p$. $(1.6.1)$, resp. $(1.6.2)$, $(1.6.1_m)$, $(1.6.2_m)$ imply
$$Y_{0}^{(1)}Q_0=Y_{0}, \eqno{(1.6.1_{0})}$$ $$Q_0X_{0}=X_{0}^{(1)}, \eqno{(1.6.2_{0})}$$
$$\Cal Y_{0}^{(1)}Q_0=\Cal Y_{0}, \eqno{(1.6.1m_0)} $$ $$Q_0\Cal X_{0}=\Cal X_{0}^{(1)} \eqno{(1.6.2m_0)}$$
Since $|Q_0|\ne0$, the Lang's theorem imply that solutions having $\det\ne0$ to $(1.6.1m_{0}), (1.6.2m_{0})$ exist and are unique up to multiplication by $GL_{r}(\n F_q)$ (from the left for $\Cal Y$, from the right for $\Cal X$). Hence, the set of solutions to $(1.6.1_{0}), (1.6.2_{0})$ is of dimension $r$ over $\n F_q$.
\medskip
After elimination of $r-1$ unknowns in $(1.6.1_\infty)$, $(1.6.2_\infty)$ we get affine equations $\g E_{(1.6.1)}$, $\g E_{(1.6.2)}$ of type (2.1) (see Appendix, 6.4 for details of elimination process). The fact that the set of solutions to $(1.6.1_{0}), (1.6.2_{0})$ is of dimension $r$ over $\n F_q$ implies that $r$ of the obtained $\g E_{(1.6.1)}$, $\g E_{(1.6.2)}$ (see (2.1)) are both our initial $r$. Further, proposition 2.3 implies that the sets of solutions to both $\g E_{(1.6.1)}$, $\g E_{(1.6.2)}$ over $\n F_q[[T]]$ also have dimension $r$. The same dimension have the sets of solutions to $(1.6.1_\infty)$, $(1.6.2_\infty)$.
\medskip
Finally, since solutions $\Cal X_{0}$, $\Cal Y_{0}$ to $(1.6.1m_{0}), (1.6.2m_{0})$ belong to $GL_r(\p)$, we have $\Cal Y_{0}\Cal X_{0}\in GL_r(\n F_q)$, hence the matrix of $\pi_\infty$ belongs to $GL_r(\n F_q[[T]])$.
$\square$
\medskip
The sets of solutions to $(1.6.1)$, $(1.6.2)$ are $\n F_q[T]$-submodules of the sets of solutions to $(1.6.1_\infty)$, $(1.6.2_\infty)$, and $\pi$ is a restriction of $\pi_\infty$ to these sets.
\medskip
{\bf 6.7.2.} Further, we have: $\n F_q[T]$-linearly independent sets of solutions to (1.6.1) (resp. (1.6.2)) are $\n F_q[[T]]$-linearly independent. Moreover, they are linearly independent also over $\p[[T]]$ (and hence over $\p\{T\}$).
\medskip
This can be proved by the same method as [P], Prop. 3.3.9(c), see Remark 2.3.5.
\medskip
{\bf Corollary 6.7.3.} In appropriate bases of the sets of solutions to $(1.6.1_\infty)$, $(1.6.2_\infty)$ the matrix of $\pi$ is a submatrix of the matrix of $\pi_\infty$.
\medskip
Formula (6.6.1) follows immediately from the fact that if we have a $c_1\times c_1$ invertible matrix $\Cal M$ and its $c_2\times c_3$ submatrix $\Cal S$ then the rank of $\Cal S$ is $\ge c_2+c_3-c_1$.
\medskip

\medskip
\medskip

{\bf References}
\medskip
[AB] S. Abhyankar, Two notes on formal power series, Proc. Amer. Math. Soc. 7 (1956), 903 -- 905.
\medskip
[AN] Anderson Greg W. $t$-motives. Duke Math. J. Volume 53, Number 2 (1986), 457 -- 502.
\medskip
[AT]  Anderson, Greg W.; Thakur, Dinesh S. Tensor powers of the Carlitz module and zeta values. Ann. of Math. (2) 132 (1990), no. 1, 159 -- 191.
\medskip
[C] C. Chevalley, Introduction to the Theory of Algebraic Functions of One Variable, Amer.Math. Soc., 1951.
\medskip
[EGL] S. Ehbauer, A. Grishkov, D. Logachev. Calculation of $h^1$ of some Anderson t-motives. To appear in Journal of Algebra and its applications. https://arxiv.org/pdf/2006.00316.pdf
\medskip
[G] Goss, D. Basic structures of function field arithmetic. Springer-Verlag, Berlin, 1996. xiv+422 pp.
\medskip
[G1] Goss, D. Private letter to D. Logachev, 2008
\medskip
[GL07] Grishkov, A., Logachev, D. Duality of Anderson t-motives. 2007. arxiv.org/pdf/0711.1928.pdf
\medskip
[GL09] Grishkov, A., Logachev, D. Anderson T-motives and abelian varieties with MIQF: results coming from an analogy. 2009.

https://arxiv.org/pdf/0907.4712.pdf
\medskip
[GL17] Grishkov, A., Logachev, D. Lattice map for Anderson t-motives: first approach. J. of Number Theory. 2017, vol. 180, p. 373 -- 402.
https://arxiv.org/pdf/1109.0679.pdf
\medskip
[HJ] U. Hartl and A.-K. Juschka, Pink's theory of Hodge structures and the Hodge conjecture
over function fields, preprint (2016), https://arxiv.org/abs/1607.01412.
\medskip
[M] A. Maurischat, Periods of t-modules as special values, preprint (2018), https://arxiv.org/abs/1802.03233
\medskip
[P] M. A. Papanikolas, Tannakian duality for Anderson-Drinfeld motives and algebraic independence of Carlitz logarithms, Invent. Math. 171 (2008), no. 1, 123 -- 174.
\enddocument